\documentclass[11pt]{article}

\usepackage{biometrics}
\usepackage{amsmath,epsfig,epsf,psfrag,bm,geometry,color}
\usepackage{ifthen,latexsym,graphics,amssymb,fancybox,fancyhdr}
\usepackage{booktabs}
\usepackage{graphicx}
\usepackage{mathrsfs}
\usepackage{rotating}
\usepackage{longtable}
\begin{document}

\title{A Note on Spatial-Temporal Lattice Modeling and  Maximum Likelihood Estimation}
\author{Xiang Zhang and Yanbing Zheng\footnote{Department of Statistics, University of Kentucky}}
\date{\today}

\maketitle
\thispagestyle{empty}

\begin{center}
\section*{Abstract}
\end{center}

Spatial-temporal linear model and the corresponding likelihood-based statistical inference are important tools for the analysis of spatial-temporal lattice data. In this paper, we study the asymptotic properties of maximum likelihood estimates under a general asymptotic framework for spatial-temporal linear models. We propose mild regularity conditions on the spatial-temporal weight matrices and derive the asymptotic properties (consistency and asymptotic normality) of maximum likelihood estimates. A simulation study is conducted to examine the finite-sample properties of the maximum likelihood estimates.

\noindent
{\it Keywords: } Autoregressive models; increasing domain asymptotics; infill asymptotics; linear regression; spatial-temporal process

\section{Introduction}

Spatial-temporal linear models are important tools for the analysis of spatial-temporal lattice data and have been applied in a wide range of disciplines (see,  e.g., \cite*{anselin01}, \cite*{baltagi05}, and \cite*{cressie93}).  A spatial-temporal linear model relates the response variable of interest to covariates via a linear regression component and models the spatial-temporal dependence in data via a random error component that is assumed to be a zero-mean Gaussian process.  In this paper, we focus on simultaneous-autoregressive (SAR) type of spatial-temporal linear models for random error components and study the asymptotic properties of statistical inference via a maximum likelihood method.

For statistical inference of spatial-temporal linear models, maximum
likelihood estimation is often adopted. For the spatial-only case,
\cite{mardiam84} established that the maximum likelihood estimates
(MLE) of the parameters are consistent and asymptotically normal as
the sample size tends to infinity for general spatial linear models.
\cite{lee04} studied the asymptotic properties of quasi-maximum
likelihood estimators (QMLE) for a lag SAR model. Recently,
\cite{zhengz11} considered an error SAR model with general
neighborhood structures and explored the asymptotic properties of
MLEs under a unified asymptotic framework. \cite{RT12} developed the asymptotic normality of
ordinary least squares (LS) and instrumental variables (IV) estimates
of linear and semiparametric partly linear regression models for
spatial data and discussed
the consistency of the estimates of the spatial covariance matrix. However, it is not clear how general the asymptotic framework is. For the
spatial-temporal case, \cite{yuj08} investigated the asymptotic
properties of QMLEs for spatial dynamic panel data with fixed
effects and proposed a bias-adjusted estimator. \cite{leey10a} and
\cite{leey10b} generalized the spatial dynamic panel data model to
include both time and individual fixed effects and studied the
asymptotic properties of QMLEs.  Here, we consider the asymptotic
properties of maximum likelihood estimation for spatial-temporal
linear models. Across space, we consider the three types of
asymptotic frameworks defined in \cite{zhengz11}, namely increasing
domain, infill and hybrid of increasing domain and infill
asymptotics.  Over time, we assume that the number of time points
tends to infinite as is traditionally done in time series, but will
discuss the case with fixed time points.

The remainder of the paper is organized as follows. In Section 2, we describe the spatial-temporal linear models. In Section 3, we develop maximum likelihood for inference and the main theoretical results are given in Section 4. A simulation study is performed in Section 5. Conclusions and discussions are given in Section 6. Technical proofs are shown in the Appendices A-B.

\section{Spatial-Temporal Linear Model}

Let $Y_{it}$ denote the response variable at site $i$ and time $t$, where $i=1,\ldots,n$ and $t=1,\ldots,m$. Let
\begin{eqnarray}
Y_{it} = \bm{X}_{it}'\bm{\beta} + \epsilon_{it},
\label{eqn:emodel}
\end{eqnarray}
where $\bm{X}_{it} = (X_{1it},\ldots,X_{pit})'$ is a vector of the covariates and $\bm{\beta}=(\beta_1,\ldots,\beta_p)'$ a vector of the regression coefficients. To formulate a spatial-temporal model for the random errors $\{\epsilon_{it}: i=1,\ldots,n, t=1,\ldots,m \}$, we focus on a SAR-type model and specify the model for $\bm{\epsilon}_t = (\epsilon_{1t},\ldots,\epsilon_{nt})'$ in terms of the errors from time $\max \{1, t-s \}$ to $t$,
\begin{eqnarray}
\bm{\epsilon}_t = \sum_{l=0}^{\min \{ s, t-1\}} \bm{C}_l \bm{\epsilon}_{t-l} + \bm{\nu}_t,
\label{eqn:sarep}
\end{eqnarray}
where $\bm{C}_l = [c^{(l)}_{ij}]_{i,j=1}^{n}$,
$l=0,\ldots,\min\{s,t-1\}$, $s \geq 0$, are spatial-temporal
dependence matrices and $\bm{\nu}_t = (\nu_{1t},\ldots,\nu_{nt})'
\sim N(0,\sigma^2 \bm{I}_n)$ is a vector of white
noises. Let $\bm{Y}_t = (Y_{1t},\ldots,Y_{nt})'$ denote the
$n-$dimensional vector of response variables  on the entire spatial
lattice for a given time point $t$ and $\bm{Y}_{nm} = (\bm{Y}_1',
\ldots,\bm{Y}_m')'$ the $nm-$dimensional vector of response
variables at $m$ time points. Let ${\bm
\epsilon}_{nm}=(\epsilon_{11},\dots,\epsilon_{n1},\dots,\epsilon_{1m},\dots,\epsilon_{nm})'$
be a $nm$-dimensional vector of random errors and ${\bm
\nu}_{nm}=(\nu_{11},\dots,\nu_{n1},\dots,\nu_{1m},\dots,\nu_{nm})'$
a $nm$-dimensional vector of white noises. We have, under
(\ref{eqn:emodel}) and (\ref{eqn:sarep}),
\[
\bm{Y}_{nm} = \bm{X}_{nm} \bm{\beta} + \bm{\epsilon}_{nm}, \quad {\rm and} \quad \bm{\epsilon}_{nm} = \bm{C} \bm{\epsilon}_{nm} + \bm{\nu}_{nm}, \]
where ${\bf X}_{nm}$ is a $nm \times p$ design matrix,  $\bm{C}$ is a lower-triangular matrix with $\bm{C}_0$ as the diagonal blocks and $\bm{C}_l$ as the $l$th sub-diagonal blocks for $l=1,\ldots,s$. So the joint distribution of the response variable is
\begin{eqnarray}
\bm{Y}_{nm} &\sim& N \left(\bm{X}_{nm}\bm{\beta}, \sigma^2 (\bm{I}_{nm} -\bm{C})^{-1} (\bm{I}_{nm} - \bm{C}')^{-1} \right).
\label{eqn:modeld}
\end{eqnarray}
where $\bm{I}_{nm}$ is an $nm \times nm$ identity matrix.

The model (\ref{eqn:sarep}) is quite general and features a variety of neighborhood structures over space and time. Let ${\cal N}_n(i)=\{j: \mbox{site $j$ is a neighbor of site $i$}\}$ denote the neighborhood of site $i$.  The neighborhood of site $i$ can be further partitioned into $q$ orders, such that ${\cal N}_n(i)=\cup_{k=1}^q{\cal N}_{n,k}(i)$ where ${\cal N}_{n,k}(i)=\{j: \mbox{site $j$ is a $k$th order neighbor of site $i$}\}$. Let ${\bf W}_{nk}=[w_{nk}^{i,j}]_{i,j=1}^{n}$ be a $n\times n$ spatial weight matrix with zero diagonal elements for all $1\leq k\leq q$. An example is binary spatial weights such that $w_{nk}^{i,j}=1$ if $j \in {\cal N}_{n,k}(i)$ and $0$ otherwise. Some special cases of model (\ref{eqn:sarep}) are as follows:
\begin{itemize}
\item Spatial independence: $s\geq 1$, $\bm{C}_0=\bm{0}$ and $\bm{C}_l = \alpha_l \bm{I}_n$ for $l=1,\ldots,s$.
\item Temporal independence: $s=0$ and $\bm{C}_0 = \sum_{k=1}^q \theta_k \bm{W}_{nk}$.
\item Spatial-temporal separable neighborhood structure: $s \geq 1, \bm{C}_0 = \sum_{k=1}^q \theta_k \bm{W}_{nk}$ and $\bm{C}_l = \alpha_l \bm{I}_{n}$ for $l=1,\ldots,s$.
\item Spatial-temporal non-separable neighborhood structure: $s \geq 1, \bm{C}_0 = \sum_{k=1}^q \theta_k \bm{W}_{nk}$ and $\bm{C}_l = \alpha_l \bm{I}_{n} + \sum_{k=1}^q \theta_k^l \bm{W}_{nk}^l$ for $l=1,\ldots,s$.
\end{itemize}
Regularity conditions on the parameter space are needed to ensure
that the model specified under (\ref{eqn:emodel}) and
(\ref{eqn:sarep}) is valid. Since $\bm{C}$ is lower-triangular, a sufficient condition to ensure the non-singularity of $\bm{I}_{nm} - \bm{C}$ is that its diagonal blocks $\bm{I}_n - \sum_{k=1}^q \theta_k \bm{W}_{nk}$  are non-singular, where $\bm{I}_{n}$ is an $n \times n$ identity matrix. For row standardized spatial weight matrices $\bm{W}_{nk}$, if $\sum_{k=1}^q |\theta_k| < 1$, then $\bm{I}_n - \sum_{k=1}^q \theta_k \bm{W}_{nk}$ is nonsingular and thus the covariance matrix in (\ref{eqn:modeld}) is positive definite (Corollary 5.6.16, \cite*{horn85}). We let $\bm{\theta}= (\theta_1,\ldots,\theta_q)'
\in \Theta$, where $\Theta$ is a compact subset of $\mathbb{R}^{q}$.
In the following, we will focus on a spatial-temporal linear model
with a separable spatial-temporal neighborhood structure and $s=1$.
The results can be extended to general $s \geq 1$ readily. We let
$\alpha \in A_{\alpha}$, where $A_{\alpha}$ is a compact set of
$(-1,1)$.

\section{Maximum Likelihood Estimations}

Let ${\bm \eta} = ({\bm \beta}',{\bm \xi}',{\sigma}^2)'$ denote the $\{p+(q+1)+1\}$-dimensional vector of unknown parameters under the model specified in (\ref{eqn:emodel}) and (\ref{eqn:sarep}), where $\bm{\xi}=(\bm{\theta}', \alpha)'$. The log-likelihood function, up to a constant, is
\begin{eqnarray}
\ell({\bm \eta})=-(nm/2)\log{\sigma}^2+\log|{\bm S}_{nm}({\bm \xi})|-(2\sigma^2)^{-1}{\bm \nu}_{nm}'{\bm \nu}_{nm}.
\label{eqn:loglike}
\end{eqnarray}
where ${\bm S}_{nm}({\bm \xi})={\bm I}_{nm}-\bm{C}$ and $\bm{\nu}_{nm} = \bm{S}_{nm}(\bm{\xi})(\bm{Y}_{nm}-\bm{X}_{nm}\bm{\beta})$.  The fist-order derivatives of $\ell({\bm \eta})$ with respect to $\bm{\beta}$ and $\sigma^2$ are, respectively,
\[ \frac{\partial \ell({\bm \eta})}{\partial {\bm \beta}} = (\sigma^2)^{-1}{\bf X}_{nm}'{\bf S}'_{nm}({\bm \xi}){\bm
\nu}_{nm}, \quad
\frac{\partial \ell({\bm \eta})}{\partial \sigma^2} = (2\sigma^4)^{-1}({\bm \nu}_{nm}'{\bm \nu}_{nm}-nm \sigma^2). \]
By setting the score functions equal to zero, we obtain the maximum likelihood estimate (MLE) of ${\bm \beta}$ and $\sigma^2$,
\begin{eqnarray*}
\hat{{\bm \beta}}_{nm}(\bm \xi) &=& \{{\bm X}_{nm}'{\bm S}'_{nm}({\bm \xi}){\bm S}_{nm}({\bm \xi}){\bm X}_{nm}\}^{-1}{\bm X}_{nm}'{\bm S}'_{nm}({\bm \xi}){\bm S}_{nm}({\bm \xi}){\bm Y}_{nm}, \\
\hat{\sigma}^2_{nm}(\bm \xi) &=& (nm)^{-1}\{ {\bm Y}_{nm}-{\bf X}_{nm}\hat{{\bm \beta}}_{nm}(\bm \xi)\}'{\bm S}'_{nm}({\bm
\xi}){\bm S}_{nm}({\bm \xi})\{ {\bm Y}_{nm}-{\bm X}_{nm}\hat{{\bm \beta}}_{nm}(\bm \xi) \}.
\end{eqnarray*}
We define a profile log-likelihood function of $\bm{\xi}$ as
\begin{eqnarray}
\ell({\bm \xi})= \ell\{\hat{{\bm \beta}}_{nm}(\bm \xi),\bm{\xi},\hat{\sigma}^2_{nm}(\bm \xi)\} = -(nm/2) \log \hat{\sigma}_{nm}(\bm{\xi})^2+\log|{\bm S}_{nm}({\bm \xi})|-nm/2.
\label{eqn:profile}
\end{eqnarray}
Then the MLE of $\bm{\xi}$ maximizes the profile log-likelihood $\ell({\bm \xi})$ and is denotes as $\hat{\bm{\xi}}_{nm}$.

\section{Asymptotic Properties}

In \cite{zhengz11}, three types of asymptotic frameworks are defined for spatial linear models on a lattice in terms of the volume of the spatial domain (i.e., the Lebesgue measure) and that of the individual cells.
\begin{itemize}
\item Increasing domain asymptotics: The volume of the spatial lattice tends to infinity while the volume of each cell on the lattice is fixed.
\item Infill asymptotics: The volume of the spatial lattice is fixed while the volume of each cell on the lattice tends to zero.
\item Hybrid asymptotics (increasing domain combined with infill asymptotics): The volume of the spatial lattice tends to infinity and the volume of each cell on the lattice tends to zero.
\end{itemize}
In this section, we study the asymptotic properties of the MLEs of the model parameters for the spatial-temporal linear model defined in (\ref{eqn:emodel}) and (\ref{eqn:sarep}). We consider all three types asymptotics in the spatial domain and assume that time tends to infinity.

Let ${\bm \eta_0} = ({\bm \beta_0}',{\bm \xi_0}',{\sigma_0}^2)'$
denote the $\{p+(q+1)+1\}$-dimensional vector of true parameters and
$\bm{S}_{0nm}=\bm{S}_{nm}(\bm{\xi}_0)$. The model evaluated at the
true parameters $\bm{\eta}_0$ is $\bm{Y}_{nm} =
\bm{X}_{nm}\bm{\beta}_0+\bm{S}_{0nm}^{-1}\bm{\nu}_{onm}$. To
establish the asymptotic properties of MLEs of the model parameters,
we impose the following regularities.
\begin{itemize}
\item[(A.1)] The elements $w^{i,j}_{nk}$ of the spatial weight matrix ${\bf W}_{nk}$ are at most of order ${h_n}^{-1}$ uniformly for all $j \neq i$ and $w^{i,j}_{nk} = 0$, for all $i = 1,\dots, n$ and $k = 1,\dots,q$ and $h_n$ is bounded away from zero uniformly.
\item[(A.2)] The sequence of spatial weight matrices $\{{\bf W}_{nk} : k = 1,\dots,q\}$ are uniformly bounded in matrix norms ${||\cdot||}_1$ and ${||\cdot||}_{\infty}$. For a matrix $\bm{A}=[a_{i,j}]^{n}_{i,j=1}$, the maximum column sum matrix norm $||\cdot||_1$ is defined by $||\bm{A}||_1 = \max_{1\leq j \leq n} \sum_{i=1}^n |a_{i,j}|$ and the maximum row sum matrix norm $||\cdot||_{\infty}$ is defined by $||\bm{A}||_{\infty} = \max_{1\leq i \leq n} \sum_{j=1}^n |a_{i,j}|$ (\cite*{horn85}).
\item[(A.3)] The matrix ${\bf S}_{n}({\bm \theta})$ is nonsingular for ${\bm \theta}\in \Theta$ and $n \in \mathbb{N}$, where ${\bf S}_{n}({\bm \theta})=\bm{I}_n - \sum_{k=1}^q \theta_k \bm{W}_{nk}$.
\item[(A.4)] The sequence of matrices $\{{\bf S}^{-1}_{n}({\bm \theta})\}$ is uniformly bounded in  matrix norms ${||\cdot||}_1$ and ${||\cdot||}_{\infty}$ for ${\bm \theta}\in \Theta$. The true parameter ${\bm \xi}_0$ is in the interior of $\Xi$, where $\Xi = \Theta \times A_{\alpha}$.
\item[(A.5)] The elements of ${\bf X}_{nm}$ are uniformly bounded constants. The limit of $({nm})^{-1}{\bf X}'_{nm}{\bf S}'_{nm}({\bm \xi}){\bf S}_{nm}({\bm \xi}){\bf X}_{nm}$ as $nm \rightarrow \infty$ exists and is nonsingular for ${\bm \xi} \in \Xi$.
\item[(A.6)] For ${\bm \xi} \neq {\bm \xi}_0$, $\lim_{nm \rightarrow \infty}h_n({nm})^{-1}\{\log|\sigma_{nm}^{*2}({\bm \xi}){\bf S}'_{nm}({\bm \xi}){\bf S}_{nm}({\bm \xi})| - \log|{\sigma}^2_{0}{\bf S}'_{0nm}{\bf S}_{0nm}|\} \neq 0$, where $\sigma_{nm}^{*2}({\bm \xi}) = (nm)^{-1}\sigma_0^2 {\rm tr}\left\{\bm{S}_{0nm}^{'-1}\bm{S}'_{nm}(\bm{\xi})\bm{S}_{nm}(\bm{\xi})\bm{S}_{0nm}^{-1} \right\}$.
\end{itemize}

Assumptions~(A.1) and (A.2) are regularity conditions on the spatial
weight matrices, where assumption~(A.2) is generally satisfied under
row standardization. Here the order of the elements in the spatial
weight matrices $h_n^{-1}$ is an essential element in the
specification of an asymptotic type in the spatial domain.  We
assume $w^{i,j}_{nk}=\mathcal{O}(h_n^{-1})$, where $h_n$ could be bounded or tend to infinity.
Consider an example with distance-based neighbors on a regular
spatial lattice (Example 2 in \cite*{zhengz11}). Let $d_{ij}$ denote
the Euclidean distance between sites $i$ and $j$ and $a^{i,j}_{nk} =
\mathcal{I}\{d_{ij} \sim (\delta_{k-1},\delta_{k}]\}$ with
prespecified threshold values $\delta_0 = 0 < \delta_1 < \dots <
\delta_q$. Then row standardized weight matrices based on
$\bm{A}_{nk}=[a_{nk}^{i,j}]_{i,j=1}^n$ are ${\bm W}_{nk}$ with
elements $w_{nk}^{i,j}=a_{nk}^{i,j}/\sum_{j=1}^n a_{nk}^{i,j}$. For
this example, $h_n = \mathcal{O}(\max\{\sum^n_{j=1} a^{i,j}_{nk} : k =
1,\dots,q, i = 1,\dots,n\})$, which is bounded under increasing
domain asymptotics. Under infill asymptotics $h_n \rightarrow
\infty$ and $h_n/n$ does not tend to 0 as $n \rightarrow \infty$,
whereas under hybrid asymptotics, $h_n \rightarrow \infty$ and
$h_n/n \rightarrow 0$ as $n \rightarrow \infty$. Assumptions~(A.3) and (A.4) are standard assumptions made
about the sequence of matrices $\bm{S}_n(\bm{\theta})$, where assumption~(A.4) is needed to ensure that the variance
of ${\bf Y}_{n}$ is bounded. For row
standardized spatial weight matrices $\bm{W}_{nk}$, if $\sum_{k=1}^q |\theta_k| < 1$, then assumption~(A.3) and (A.4) are satisfied (Corollary 5.6.16, \cite*{horn85}). Assumption~(A.5) is a standard assumption of the design matrix and
implies that the elements of $nm \{\bm{X}'_{nm}
\bm{S}'_{nm}(\bm{\xi}) \bm{S}_{nm}(\bm{\xi}) \bm{X}_{nm}\}^{-1}$ are
uniformly bounded. Assumption~(A.6) is needed to establish identifiable uniqueness when
establishing consistency of the MLE (see, e.g. \cite*{lee04}, \cite*{zhengz11}). For the case with spatial independence, it can be shown that assumption~(A.6) is simplified to $\lim_{nm \rightarrow \infty}h_n \log \left[(nm)^{-1}{\rm tr}\left\{\bm{S}_{0nm}^{'-1}\bm{S}'_{nm}(\bm{\xi})\bm{S}_{nm}(\bm{\xi})\bm{S}_{0nm}^{-1} \right\}\right] \neq 0$ for ${\bm \xi} \neq {\bm \xi}_0$, which is satisfied if $h_n$ converges to a non-zero point. 

Let $\hat{{\bm \beta}} = \hat{{\bm
\beta}}_{nm}(\hat{\bm{\xi}}_{nm}),  \hat{\sigma}_{nm}^2 =
\hat{\sigma}_{nm}^2(\hat{\bm{\xi}}_{nm})$, and $\hat{{\bm
\eta}}_{nm} = (\hat{{\bm \beta}}_{nm}',\hat{{\bm
\xi}}_{nm}',\hat{\sigma}_{nm}^2)'$ be the MLE of ${\bm \eta}$. We
have the asymptotic properties of the MLE of ${\bm \eta}$ as
follows.
\begin{Thm}
Assume that (A.1)-(A.6) hold and $h_n/(nm) \rightarrow 0$ as $nm \rightarrow \infty$. Then the MLE of ${\bm \eta}$ is consistent such that, as $nm \rightarrow \infty$, $\hat{{\bm \eta}}_{nm} \stackrel{\mathcal{P}}{\rightarrow} \bm{\eta}_0$.
\label{thm:consistency}
\end{Thm}
Theorem \ref{thm:consistency} gives the consistency of the MLE of model parameters. It shows that, under the
regularity conditions, when the size of the spatial domain $n \rightarrow \infty$,
\begin{itemize}
\item
if the number of time points $m$ is fixed, then the MLE of
$\bm{\eta}$ is consistent when either $h_n=\mathcal{O}(1)$ or $h_n
\rightarrow \infty$ but $h_n/n \rightarrow 0$ as $n \rightarrow
\infty$, which correspond to increasing domain asymptotics and
hybrid asymptotics in the spatial domain, respectively. This result
is the same as for the spatial-only case in \cite{zhengz11}.
\item
if the number of time points $m \rightarrow \infty$, the MLE of
$\bm{\eta}$ is consistent in cases (i) $h_n=\mathcal{O}(1)$; (ii)
$h_n \rightarrow \infty$ but $h_n/n \rightarrow 0$; or (iii)
$h_n\rightarrow \infty$ but $h_n/n \rightarrow c\in (0,\infty]$,
where the case (iii) corresponds to infill asymptotics in the
spatial domain.
\end{itemize}
When the size of spatial domain $n$ is fixed, the MLE of
$\bm{\eta}$ is consistent only if $m \rightarrow \infty$, which is consistent with results in \cite{yuj08}.

\begin{Thm}
Assume that (A.1)-(A.6) hold.
\begin{itemize}
\item[(i)] If $h_n = {\cal O}(1)$ and the limit of $-({nm})^{-1}E \left\{\frac{\partial^2 \ell({\bm \eta})}{\partial{{\bm \eta}}\partial{\bm \eta}'}\right\}$ as $nm \rightarrow \infty$ exists and is positive definite for $\bm{\eta} \in \mathbb{R}^{p}\times \Xi \times \mathbb{R}^{+}$, then the MLE of ${\bm \eta}$ is asymptotic normal such that, as $nm \rightarrow \infty$,
\[({nm})^{1/2} (\hat{{\bm \eta}}_{nm} - {\bm \eta}_0) \stackrel{\mathcal{D}}{\rightarrow} N(\bm{0},{\bm \Sigma}_{\bm{\eta}_0})\]
where ${\bm \Sigma}^{-1}_{\bm{\eta}_0} =-\lim_{nm \rightarrow \infty}({nm})^{-1}E
\left\{\frac{\partial^2 \ell({\bm \eta}_0)}{\partial{{\bm
\eta}}\partial{\bm \eta}'}\right\}$.
\item[(ii)] If $h_n \rightarrow \infty$ and $h_n^{1+\delta}/(nm) \rightarrow 0$ for some $\delta > 0$ as $nm \rightarrow \infty$, and if the limit of $-h_n({nm})^{-1}E \left\{\frac{\partial^2 \ell({\bm \eta})}{\partial{{\bm \theta}}\partial{\bm \theta}'}\right\}$ as $nm \rightarrow \infty$ exists and is positive definite for $\bm{\theta} \in \Theta$, then the MLE of $\bm{\eta}$ is asymptotically normal such that, as $nm \rightarrow \infty$,
\begin{eqnarray*}
(nm)^{1/2}(\hat{\bm{\beta}}_{nm} - \bm{\beta}_0) &\stackrel{\mathcal{D}}{\rightarrow}& N \left(\bm{0}, \bm{\Sigma}_{\bm{\beta}_0}\right),
\quad (nm)^{1/2}(\hat{\sigma}^2_{nm} - \sigma^2_0) \stackrel{\mathcal{D}}{\rightarrow} N(0, 2\sigma^4_0) \\
(nm/h_n)^{1/2}(\hat{\bm{\theta}}_{nm} - \bm{\theta}_0) &\stackrel{\mathcal{D}}{\rightarrow}& N(\bm{0}, \bm{\Sigma}_{\bm{\theta}_0}),
\quad (nm)^{1/2}(\hat{\alpha}^2_{nm} - \alpha^2_0) \stackrel{\mathcal{D}}{\rightarrow} N\left(0, \Sigma_{\alpha_0}\right),
\end{eqnarray*}
where  $\bm{\Sigma}_{\bm{\beta}_0} = \sigma^2_0 \lim_{nm \rightarrow \infty} nm\left(\bm{X}_{nm}' \bm{S}_{0nm}' \bm{S}_{0nm} \bm{X}_{nm} \right)^{-1}$, $\bm{\Sigma}_{\bm{\theta}_0}^{-1} = \lim_{nm \rightarrow \infty} E \left\{-h_n(nm)^{-1}\frac{\partial^2 \ell(\bm{\eta}_0)}{\partial \bm{\theta} \partial \bm{\theta}'}  \right\}$ and $\Sigma_{\alpha_0} = \lim_{nm \rightarrow \infty} nm \left\{{\rm tr}(\bm{H}_{nm}^2+\bm{H}_{nm}'\bm{H}_{nm}) \right\}^{-1}$. Here $\bm{H}_{nm}=\bm{F}_{nm} \bm{S}_{0nm}^{-1}$, where $\bm{F}_{nm}$ is a matrix with all elements equal to zero except that the first-order lower sub-diagonal blocks are $\bm{I}_n$.
\end{itemize}
\label{thm:normality}
\end{Thm}
Theorem \ref{thm:normality} is about the normality of the MLE of model parameters. Theorem \ref{thm:normality} (i) shows that
\begin{itemize}
\item
if the sample size $n \rightarrow \infty$ and time points $m
\rightarrow \infty$, then the MLE of $\bm{\eta}$ is asymptotically normal
at a convergence rate of square root of $nm$ when $h_n =
\mathcal{O}(1)$.
\item
if the time points $m$ is fixed but $n \rightarrow \infty$, the
convergence rate is reduced to the square root of $n$ as for the
spatial-only case.
\item
if the size of spatial lattice $n$ is fixed but $m \rightarrow
\infty$, then the MLE of $\bm{\eta}$ is asymptotically normal at a
convergence rate of square root of $m$, which is consistent with results in \cite{yuj08}.
\end{itemize}
Theorem \ref{thm:normality} (ii) requires that the size of spatial lattice $n \rightarrow \infty$. It shows that
\begin{itemize}
\item
if $m \rightarrow \infty$, when $h_n \rightarrow \infty$ but $h_n/n
\rightarrow 0$, or $h_n\rightarrow \infty$ but $h_n/n \rightarrow
c\in (0,\infty]$, then the MLE of $\bm{\eta}$ is asymptotically
normal with a convergence rate $\sqrt{nm}$ for the regression
coefficients $\bm{\beta}$, the temporal autoregressive coefficient
$\alpha$ and the variance component $\sigma^2$ and a convergence
rate $\sqrt{nm/h_n}$ for the spatial autoregressive coefficients
$\bm{\theta}$.
\item
if $m$ is fixed, then similarly to the spatial-only case in
\cite{zhengz11}, the MLE
of $\bm{\eta}$ is asymptotically normal with a reduced convergence
rate $\sqrt{n}$ for the regression coefficients $\bm{\beta}$, the
temporal autoregressive coefficient $\alpha$ and the variance
component $\sigma^2$ and a convergence rate $\sqrt{n/h_n}$ for the
spatial autoregressive coefficients $\bm{\theta}$ when $h_n
\rightarrow \infty$ and $h_n/n \rightarrow 0$.
\end{itemize}

The MLEs of model parameters can be obtained using the Newton-Raphson algorithm. For the large-sample case, according to the consistency of MLEs, under the assumptions in Theorem 2 (i) we can estimate the covariance matrix $\bm{\Sigma}_{\bm{\eta}_0}$ by using
\[ \left[-\left.({nm})^{-1}E
\left\{\frac{\partial^2 \ell({\bm \eta})}{\partial{{\bm
\eta}}\partial{\bm \eta}'}\right\}\right]^{-1}\right|_{\bm{\eta}=\hat{\bm{\eta}}_{nm}}. \]
The detailed formats of the elements in $-({nm})^{-1}E
\left\{\frac{\partial^2 \ell({\bm \eta})}{\partial{{\bm
\eta}}\partial{\bm \eta}'}\right\}$ are given in Appendix B. The covariance structure under the assumptions in Theorem 2 (ii) can be estimated similarly. 

\section{Simulation}

We now conduct a simulation study to examine the finite-sample properties of the MLEs under the three types of asymptotics across space and with an increasing number of time points. We consider an $r \times r$ square lattice with a unit resolution and $m$ temporal points. We vary the number of time points by letting $m=2$, $5$ {\bf or $10$}. For each value of $m$, we vary the lattice size by letting $r = 4$ or $8$. For each lattice size, we further divide each cell into an $r^* \times r^*$ sub-lattice and vary the sub-lattice size by letting $r^*= 1, 2$, or $4$. Thus, for each time point, the sample size $n$ ranges from $16$ ($r = 4$; $r^* =1$) to $1024$ ($r = 8$; $r^*= 4$).

For a given lattice size $r$, sub-lattice size $r^*$ and temporal
length $m$, we simulate data from the spatial-temporal model defined
in (\ref{eqn:emodel}) and (\ref{eqn:sarep}).  For the linear
regression, we let $E(Y_{it})=\beta_0+\beta_1 X_i$, where
$X_i=\sin(i)$, $\beta_0=2$, and $\beta_1=2$ for the $i$th cell,
$i=1,\ldots,n$ and $t=1,\ldots,m$.  For the spatial dependence, we
consider distance-based neighborhood with order $q=1$.  We let
$a_{n1}^{i,j}=\mathcal{I}\{d_{ij} \in (0,1]\}$, where $d_{ij}$
denotes the Euclidean distance between sites $i$ and $j$, and then
define a row standardized weight matrix $w_{n1}^{i,j}=a_{n1}^{i,j}
/\sum_{j=1}^n a_{n1}^{i,j}$. The parameter values are set at
$\theta_1=0.8$, $\alpha=0.2$, and $\sigma^2=1$.  For each simulated
data, we estimate the model parameters by using the maximum likelihood method and
obtain $\hat{\beta}_0, \hat{\beta}_1, \hat{\theta}_1$,
$\hat{\alpha}$, and $\hat\sigma^2$.  We repeat this procedure 100
times.

Tables \ref{tab:simulation1}, \ref{tab:simulation2} and  \ref{tab: simulation3} give the
means and standard deviations of the MLEs. First, for a given number
of time points $m$, we note that in general, the biases and standard
deviations of all five parameter estimates decrease as the lattice
size $r$ increases from 4 to 8 for any given sub-lattice size $r^*$
or as both $r$ and $r^*$ increase, which correspond to the
increasing domain asymptotics and hybrid asymptotics across,
respectively. Next, for a given number of time points $m$ and
given lattice size $r$, we consider the results over all sub-lattice
size $r^*$, which corresponds to infill asymptotics. In general, the
biases and standard deviations of the regression coefficient estimates
$\hat{\beta}_0$, $\hat{\beta_1}$, the variance component estimate
$\hat{\sigma}^2$ and the temporal autoregressive coefficient estimate
$\hat{\alpha}$ decrease as the sub-lattice size $r^*$ increases from
$1 \times 1$ to $4 \times 4$. However, for the spatial
autoregressive coefficient estimate $\hat{\theta}$, its biases and
standard deviations remain similar as $r^*$ increases, which indicates that $\hat{\theta}$ is inconsistent under the infill asymptotics when the number of 
time points $m$ is fixed. This result agrees with the spatial-only case in
\cite{zhengz11}. Last, as the number of time points $m$ increases
from 2 to 10, we note that the overall performance of all five
parameter estimates improves with either fixed or increasing lattice
size $r$ and sub-lattice size $r^*$.

\section{Conclusions and Discussion}

In this paper, we have studied the asymptotic properties of maximum likelihood estimates under a general asymptotic framework for spatial-temporal linear models. We have considered three types of asymptotics in the spatial domain and let the number of time points tend to infinity. Under mild regularity conditions on the spatial-temporal weight matrices, we have derived the asymptotic properties (consistency and asymptotic normality) of maximum likelihood estimates. The results can be easily extended to models with temporal lags $s > 1$.  It is plausible that the asymptotics of MLEs for models with a general non-separable spatial-temporal neighborhood structure can be developed in a similar technique, which is currently under investigation.

In our spatial-temporal autoregressive models, we assume that the errors are zero with $\bm{\epsilon}_l=\bm{0}$ at initial time points $1-s \leq l \leq 0$. An alternative way to formulate the process is to pre-specify a distribution for the errors at the initial time points. An analogy is an AR(1) model in time series $\epsilon_t=\rho \epsilon_{t-1}+\nu_t$, where $\nu_t \sim {\rm iid}~N(0, \sigma^2_{\nu})$. It is conventional to let $\epsilon_1 \sim N(0, \sigma^2_{\nu}/(1-\rho^2))$ such that ${\rm var}(\epsilon_t) = \sigma^2_{\nu}/(1-\rho^2)$. However, this would be challenging for general spatial-temporal process. In addition, a spatial-temporal process can be proposed to condition on the initial $s$ time points with $\bm{\epsilon}_t = \sum_{l=0}^s \bm{C}_l \bm{\epsilon}_{t-l}+\bm{\nu}_t$, where $\bm{\nu}_t=(\nu_{1t},\ldots,\nu_{nt})' \sim {\rm iid}~N(0,\sigma^2)$ for $t=s+1,\ldots,m$.  When $m$ goes to infinity, MLEs will behave similarly to those under our model specification.

\section*{Appendix A: Proof of Theorem 1}

From (\ref{eqn:loglike}), we have that, under the true parameters $\bm{\eta}_0$,
\begin{eqnarray*}
E \{ \ell(\bm{\eta}) \} &=&  - (nm/2)\log \sigma^2 + \log | \bm{S}_{nm}(\bm{\xi}) |\\
&& - (2\sigma^2)^{-1} \left[ (\bm{\beta}_0 - \bm{\beta})'\bm{X}_{nm}'\bm{S}'_{nm}(\bm{\xi})\bm{S}_{nm}(\bm{\xi}) \bm{X}_{nm}(\bm{\beta}_0 - \bm{\beta}) + \sigma^2_0 {\rm tr}\left\{\bm{S}_{0nm}^{'-1}\bm{S}'_{nm}(\bm{\xi})\bm{S}_{nm}(\bm{\xi})\bm{S}_{0nm}^{-1}  \right\}  \right].
\end{eqnarray*}
The first derivatives of $E \ell({\bm \eta})$ with respect to ${\bm \beta}$ and $\sigma^2$ are, respectively,
\begin{eqnarray*}
\frac{\partial E \ell({\bm \eta})}{\partial {\bm \beta}} = (\sigma^2)^{-1}{\bf X}_{nm}'{\bf S}'_{nm}({\bm \xi})\bm{X}_{nm}({\bm \beta_0}-{\bm \beta}), \quad
\frac{\partial E \ell({\bm \eta})}{\partial \sigma^2} = (2\sigma^4)^{-1}\{ E({\bm \nu}_{nm}'{\bm \nu}_{nm})-nm \sigma^2 \}.
\end{eqnarray*}
Thus the maximizers of $E \ell({\bm \eta})$ are
\[ {\bm \beta}^*_{nm}(\bm{\xi}) = {\bm \beta}_0, \quad {\rm and} \quad \sigma_{nm}^{*2}(\bm{\xi}) = (nm)^{-1} {\sigma_0}^2 {\rm tr}\{ {{{\bm S}^{'-1}_{0nm}}}{\bf S}'_{nm}({\bm \xi}){\bf S}_{nm}({\bm \xi}){\bf S}_{0nm}^{-1}\}. \]
Let $g_{nm}({\bm \xi}) = E \left[\ell\{{\bm \beta}^*_{nm}(\bm{\xi}),{\bm \xi},\sigma_{nm}^{*2}(\bm{\xi}) \}\right]=-(nm/2) \log \sigma_{nm}^{*2}({\bm \xi}) + \log|{\bm S_{nm}}({\bm \xi})|- nm/2$. We establish the consistency of $\hat{{\bm \xi}}_{nm}$ by showing that $\underset {{\bm \xi}\in \Xi }{\sup}~h_n(nm)^{-1}|\ell({\bm \xi})-g_{nm}({\bm \xi})| = o_{p}(1)$ and that $h_n(nm)^{-1}g_{nm}(\bm{\xi})$ is identifiably unique (\cite*{white94}).

To show $\underset {{\bm \xi}\in \Xi }{\sup}~h_n(nm)^{-1}|\ell({\bm \xi})-g_{nm}({\bm \xi})| = o_{p}(1)$, we have
\begin{eqnarray*}
h_n(nm)^{-1}(\ell({\bm \xi})-g_{nm}({\bm \xi}))
=-\frac{h_n}{2}\left\{\log{\hat{\sigma}_{nm}}^2({\bm \xi}) -
\log \sigma_{nm}^{*2}(\bm{\xi})
\right\}=-h_n\{2\tilde{\sigma}_{nm}^2({\bm
\xi})\}^{-1}\left\{{\hat{\sigma}_{nm}}^2({\bm \xi})-
\sigma_{nm}^{*2}({\bm \xi})\right\},
\end{eqnarray*}
where $\tilde{\sigma}_{nm}^2({\bm \xi}) = \lambda
\sigma_{nm}^{*2}(\bm{\xi}) + (1-\lambda){\hat{\sigma}_{nm}}^2({\bm
\xi})$ for some $\lambda \in (0,1)$ and $\hat{\sigma}_{nm}^2({\bm
\xi}) = (nm)^{-1}{{\bm \nu}'_{0nm}}{\bm B}_{nm}({\bm \xi}){\bm
\nu}_{0nm}$ with
\begin{eqnarray*}
{\bm B}_{nm}({\bm \xi}) = {\bm{S}_{0nm}^{'-1}}{\bm S}'_{nm}({\bm \xi})[{\bm I}_{nm} - {\bm S}_{nm}({\bm \xi}){\bm X}_{nm}\{{{\bm X}'_{nm}}{\bm S}'_{nm}({\bm \xi}){\bm S}_{nm}({\bm \xi}){\bm X}_{nm}\}^{-1}{{\bm X}'_{nm}}{\bm S}'_{nm}({\bm \xi})]{\bm S}_{nm}({\bm \xi}){{\bm S}^{-1}_{0nm}}.
\end{eqnarray*}
First, we have that $h_n \left\{\hat{\sigma}^2_{nm}(\bm \xi) - \sigma_{nm}^{*2}({\bm \xi})\right\}=o_p(1)$. By (A.1)-(A.4), we know
\begin{eqnarray*}
h_n(nm)^{-1}\left[{{\bm \nu}'_{0nm}}{\bm{S}_{0nm}^{'-1}}{\bm S}'_{nm}({\bm \xi}){\bm S}_{nm}({\bm \xi}){{\bm S}^{-1}_{0nm}}{\bm \nu}_{0nm} - E\left\{{{\bm \nu}'_{0nm}}{\bm{S}_{0nm}^{'-1}}{\bm S}'_{nm}({\bm \xi}){\bm S}_{nm}({\bm \xi}){{\bm S}^{-1}_{0nm}}{\bm \nu}_{0nm}\right\}\right] = o_{p}(1),
\end{eqnarray*}
as $nm/h_n \rightarrow \infty$, where the convergence is uniform on $\Xi$ because of the linear-quadratic form in ${\bm \xi}$ and by corollary 2.2 of \cite{newey91}.

By (A.1)-(A.5),
\begin{eqnarray*}
&&h_n(nm)^{-1}{{\bm \nu}'_{0nm}}{\bm{S}_{0nm}^{'-1}}{\bm S}'_{nm}({\bm \xi}){\bm S}_{nm}({\bm \xi}){\bm X}_{nm}\left\{{{\bm X}'_{nm}}{\bm S}'_{nm}({\bm \xi}){\bm S}_{nm}({\bm \xi}){\bm X}_{nm}\right\}^{-1}{{\bm X}'_{nm}}{\bm S}'_{nm}({\bm \xi}){\bm S}_{nm}({\bm \xi}){{\bm S}^{-1}_{0nm}}{\bm \nu}_{0nm}\\
&=& h_n(nm)^{-1}\left\{(nm)^{-1/2}{{\bm X}'_{nm}}{\bm S}'_{nm}({\bm \xi}){\bm S}_{nm}({\bm \xi}){{\bm S}^{-1}_{0nm}}{\bm \nu}_{0nm}\right\}'\left\{(nm)^{-1}{{\bm X}'_{nm}}{\bm S}'_{nm}({\bm \xi}){\bm S}_{nm}({\bm \xi}){\bm X}_{nm}\right\}^{-1}\\
&&\times \left\{(nm)^{-1/2}{{\bm X}'_{nm}}{\bm S}'_{nm}({\bm \xi}){\bm S}_{nm}({\bm \xi}){{\bm S}^{-1}_{0nm}}{\bm \nu}_{0nm}\right\}\\
&=& o_{p}(1)
\end{eqnarray*}
Again, the convergence is uniform on $\Xi$, since $(nm)^{-1}{{\bm X}'_{nm}}{\bm S}'_{nm}({\bm \xi}){\bm S}_{nm}({\bm \xi}){{\bm S}^{-1}_{0nm}}{\bm \nu}_{0nm} = o_{p}(1)$ uniformly on $\Xi$, $\{(nm)^{-1}{{\bm X}'_{nm}}{\bm S}'_{nm}({\bm \xi}){\bm S}_{nm}({\bm \xi}){\bm X}_{nm}\}^{-1}$ is uniformly bounded in $\textit{l}_{\infty}$, and the boundedness is uniform on $\Xi$ by (A.5). Thus it follows that, uniformly on $\Xi$,
\begin{eqnarray*}
&&h_n \left\{\hat{\sigma}^2_{nm}(\bm \xi) - \sigma_{nm}^{*2}({\bm \xi})\right\}\\
&=& h_n(nm)^{-1}[{{\bm \nu}'_{0nm}}{\bm{S}_{0nm}^{'-1}}{\bm
S}'_{nm}({\bm \xi}){\bm S}_{nm}({\bm \xi}){{\bm S}^{-1}_{0nm}}{\bm
\nu}_{0nm}
- E\{{{\bm \nu}'_{0nm}}{\bm{S}_{0nm}^{'-1}}{\bm S}'_{nm}({\bm \xi}){\bm S}_{nm}({\bm \xi}){{\bm S}^{-1}_{0nm}}{\bm \nu}_{0nm}\}]\nonumber\\
&&- h_n(nm)^{-1}{{\bm \nu}'_{0nm}}{\bm{S}_{0nm}^{'-1}}{\bm S}'_{nm}({\bm \xi}){\bm S}_{nm}({\bm \xi}){\bm X}_{nm}\{{{\bm X}'_{nm}}{\bm S}'_{nm}({\bm \xi}){\bm S}_{nm}({\bm \xi}){\bm X}_{nm}\}^{-1}{{\bm X}'_{nm}}{\bm S}'_{nm}({\bm \xi}){\bm S}_{nm}({\bm \xi}){{\bm S}^{-1}_{0nm}}{\bm \nu}_{0nm}\nonumber\\
&=& o_{p}(1)
\end{eqnarray*}
Next, we show the uniform boundedness of $\tilde{\sigma}_{nm}^2({\bm \xi})$. By Jensen's inequality,
\begin{eqnarray}
(nm)^{-1}\{g_{nm}({\bm \xi})-g_{nm}({\bm \xi}_0)\} &=& (nm)^{-1}(\log|{\bm S_{nm}}(\bm \xi)|-\log|{\bm S}_{0nm}|)-1/2\{\log \sigma_{nm}^{*2}({\bm \xi}) - \log{\sigma_0}^2\} \nonumber\\
&=& \frac{1}{2}\log\frac{|{\bm S}^{'-1}_{0nm}{\bm
S}'_{nm}({\bm \xi}){\bm S}_{nm}({\bm \xi}){{\bm
S}^{-1}_{0nm}}|^{nm^{-1}}}{(nm)^{-1}{\rm tr}\{ {\bm
S}^{'-1}_{0nm}{\bm S}'_{nm}({\bm \xi}){\bm S}_{nm}({\bm \xi}){{\bm
S}^{-1}_{0nm}}\}} \leq 0 \label{eqn:gjensen}
\end{eqnarray}
for ${\bm \xi} \in \Xi$. Under (A.1)-(A.4),
\begin{eqnarray}
&&(nm)^{-1}(\log|{\bm S}_{nm}({\bm \xi})|-\log|{\bm S}_{0nm}|) \nonumber\\
&=&- (nm)^{-1}\left[{\rm tr}\{{\bm S}_{nm}^{-1}(\tilde{{\bm \xi}}){\rm diag}\{\bm{I}_m \otimes {\bm W}_{n1}\},\ldots,{\rm tr}\{{\bm S}_{nm}^{-1}(\tilde{{\bm \xi}}){\rm diag}\{\bm{I}_m \otimes {\bm W}_{nq}\},{\rm tr}\{\frac{1}{\alpha}{\bm S}_{nm}^{-1}(\tilde{{\bm \xi}}){\bm A}(\alpha)\}\right]({\bm \xi}-{\bm \xi}_0)\nonumber\\
&=&-\sum^q_{k=1}{\cal O}({h_n}^{-1})({\theta}_k - {\theta}_{0k})
\label{eqn:staylor}
\end{eqnarray}
where $\tilde{{\bm \xi}}=\lambda{\bm \xi} + (1-\lambda){\bm \xi}_0$ for some $\lambda \in (0,1)$ and $\bm{A}(\alpha)=\alpha \bm{F}_{nm}$. Thus,
\begin{eqnarray*}
\log \sigma^{*2}_{nm}({\bm \xi}) &=& -2(nm)^{-1}\{g_{nm}({\bm \xi})-g_{nm}({\bm \xi}_0)\} + 2(nm)^{-1}\{\log|{\bm S}_{nm}(\bm \xi)|-\log|{\bm S}_{0nm}|\}+ \log {\sigma_0}^2\\
& \geq & 2(nm)^{-1}\{\log|{\bm S}_{nm}(\bm \xi)|-\log|{\bm S}_{0nm}|\} + \log {\sigma_0}^2
\end{eqnarray*}
which is bounded from below uniformly on  $\Xi$, and  $\sigma_{nm}^{*2}({\bm \xi})$ is bounded away from $0$ uniformly on $\Xi$. Since $\hat{\sigma}^2_{nm}(\bm{\xi}) - \sigma_{nm}^{*2}(\bm{\xi}) = o_{p}(1)$ uniformly on $\Xi$, $\hat{\sigma}^2_{nm}(\bm{\xi})$ is bounded away from $0$ in probability uniformly on $\Xi$. Hence, $\underset {\bm{\xi}\in \Xi }{\sup}~h_n(nm)^{-1}|\ell({\bm \xi})-g_{nm}({\bm \xi})| = o_{p}(1)$.

To show the identifiable uniqueness of $h_n(nm)^{-1}g_{nm}({\bm \xi})$, we note that $h_n(nm)^{-1}g_{nm}({\bm \xi})$ is uniformly equicontinuous. In
\begin{eqnarray*}
&& h_n(nm)^{-1} \{g_{nm}(\bm{\xi}_1)-g_{nm}(\bm{\xi}_2)\}\\
&=& h_n(nm)^{-1} \{ \log |\bm{S}_{nm}(\bm{\xi}_1)| - \log |\bm{S}_{nm}(\bm{\xi}_2) \} - h_n/2 \{ \log \sigma_{nm}^{*2}(\bm{\xi}_1)- \log \sigma_{nm}^{*2}(\bm{\xi}_2)\}\\
&=& h_n(nm)^{-1}\{\log| \bm{S}_{nm}(\bm{\xi}_1)| -\log|\bm{S}_{nm}(\bm{\xi}_2) \} -h_n(2 \bar{\sigma}_{nm}^{*2})^{-1}\{\sigma_{nm}^{*2}(\bm{\xi}_1) - \sigma_{nm}^{*2}(\bm{\xi}_2) \},
\end{eqnarray*}
where $\bar{\sigma}_{nm}^{*2}= \lambda \sigma_{nm}^{*2}({\bm \xi}_1) + (1-\lambda)\sigma_{nm}^{*2}({\bm \xi}_2)$ for some $\lambda \in (0,1)$ and is bounded away from $0$, both terms are uniformly equicontinuous. Since by (\ref{eqn:staylor}), $h_n(nm)^{-1}(\log|{\bm S}_{nm}({\bm \xi}_1)| -\log|{\bm S}_{nm}({\bm \xi}_2)) = -{\cal O}(1)\sum^q_{k=1}({\theta}_k - {\theta}_{0k})$ and with $\tilde{\bm \xi} = \lambda{\bm \xi}_1 + (1-\lambda){\bm \xi}_2$ for some $\lambda \in (0,1)$, we have
\begin{eqnarray*}
&& h_n \left\{\sigma_{nm}^{*2}({\bm \xi}_1) - \sigma_{nm}^{*2}({\bm \xi}_2)\right\}\\
&=& -h_n(nm)^{-1}{\sigma_0}^2 \sum^q_{k=1}{\rm tr}\left\{ {\bm{S}_{0nm}^{'-1}}{\rm diag}(\bm{I}_m \otimes {\bm W}'_{nk}){\bm S}_{nm}(\tilde{\bm \xi}){{\bm S}^{-1}_{0nm}}\right.\\
&&\quad\quad\left.+ {\bm{S}_{0nm}^{'-1}}{\bm S}'_{nm}(\tilde{\bm \xi}){\rm diag}(\bm{I}_m \otimes {\bm W}_{nk}){{\bm S}^{-1}_{0nm}}\right\}(\theta_{1k}-\theta_{2k})+\sigma^2_0{\cal O}(1)(\alpha_1-\alpha_2)\\
&=&- {\sigma_0}^2 {\cal O}(1) \sum^q_{k=1}(\theta_{1k} -
\theta_{2k})+\sigma^2_0{\cal O}(1)(\alpha_1-\alpha_2).
\end{eqnarray*}
Thus together with (A.6) and (\ref{eqn:gjensen}), $h_n(nm)^{-1}g_{nm}({\bm \xi})$ is identifiably unique. Thus, the MLE of ${\bm \xi}$ is a consistent estimator.

The consistency of $\hat{\sigma}^2_{nm}(\hat{\bm{\xi}}_{nm})$ can be derived directly from the consistency of $\sigma_{nm}^{*2}(\hat{\bm{\xi}}_{nm})$. Further,
\begin{eqnarray*}
&&\hat{\bm{\beta}}_{nm}(\hat{\bm{\xi}}_{nm})\\
&=& {\bm \beta}_0 + \{{\bm X}'_{nm}{\bm S}'_{nm}(\hat{\bm{\xi}}_{nm}){\bm S}_{nm}(\hat{\bm{\xi}}_{nm}){\bm X}_{nm}\}^{-1}{\bm X}'_{nm}{\bm S}'_{nm}(\hat{{\bm \xi}}_{nm}){\bm S}_{nm}(\hat{\bm{\xi}}_{nm}){\bm S}^{-1}_{0nm}{\bm \nu}_{0nm}\\
&= &  {\bm \beta}_0 + \{{\bm X}'_{nm}{\bm S}'_{nm}(\hat{\bm{\xi}}_{nm}){\bm S}_{nm}(\hat{\bm{\xi}}_{nm}){\bm X}_{nm}\}^{-1}{\bm X}'_{nm}{\bm S}'_{nm}(\hat{\bm{\xi}}_{nm}){\bm \nu}_{0nm} \\
&& + \sum^q_{k=1}(\theta_{0k} - \hat{\theta}_{nmk}) \{{\bm X}'_{nm}{\bm S}'_{nm}(\hat{\bm{\xi}}_{nm}){\bm S}_{nm}(\hat{\bm{\xi}}_{nm}){\bm X}_{nm}\}^{-1}{\bm X}'_{nm}{\bm S}'_{nm}(\hat{\bm{\xi}}_{nm}){\rm diag}(\bm{I}_m \otimes {\bm W}_{nk}){\bm S}^{-1}_{0nm}{\bm \nu}_{0nm}\\
& &+  \{{\bm X}'_{nm}{\bm S}'_{nm}(\hat{\bm{\xi}}_{nm}){\bm S}_{nm}(\hat{\bm{\xi}}_{nm}){\bm X}_{nm}\}^{-1}{\bm X}'_{nm}{\rm diag}(\hat{\alpha}_{nm}^2{\bm I}_{n},\ldots,\hat{\alpha}_{nm}^2{\bm I}_{n},\mathbf{0}){\bm S}^{-1}_{0nm}{\bm \nu}_{0nm}
\end{eqnarray*}
where the last three terms are of order $o_p(1)$ by (A.1)-(A.5).

\section*{Appendix B: Proof of Theorem 2}

\subsection*{The case $h_n = {\cal O}(1)$}

By (A.4), $\bm{\xi}_0$ is in the interior of $\Xi$. Thus, for sufficiently small $\epsilon > 0$, we have $A_{\epsilon}=\{ \bm{\eta}: || \bm{\eta}-\bm{\eta}_0|| < \epsilon\} \subset \mathbb{R}^p \times \Xi \times \mathbb{R}^{+}$ and $P(\hat{\bm{\eta}}_{nm} \in A_{\epsilon}) \rightarrow 1$ as $n \rightarrow \infty$, where $||\cdot||$ denotes the Euclidean norm. We establish the asymptotic normality of the MLE by showing asymptotic normality of $(nm)^{-1}\frac{\partial \ell({\bm{\eta}_0})}{\partial \bm{\eta}}$ and convergence in probability of $(nm)^{-1}\frac{\partial^2 \ell({\tilde{\bm{\eta}}_{nm}})}{\partial \bm{\eta} \partial \bm{\eta}'}$, where $\tilde{\bm{\eta}}_{nm}=\lambda \bm{\eta}_0 + (1-\lambda)\hat{\bm{\eta}}_{nm}$ for $\lambda \in (0,1)$ converges to $\bm{\eta}_0$ in probability.

For convergence of $(nm)^{-1}\frac{\partial^2 \ell(\tilde{{\bm
\eta}}_{nm})}{\partial{{\bm \eta}}\partial{\bm \eta}'}$, we show
that, under (A.1)-(A.5), $(nm)^{-1}\{\frac{\partial^2 \ell(\tilde{\bm
\eta}_{nm})}{\partial{{\bm \eta}}\partial{\bm \eta}'} -
\frac{\partial^2 \ell({\bm \eta}_0)}{\partial{{\bm \eta}}\partial{\bm
\eta}'}\} = o_p(1)$ and ${nm}^{-1}\{\frac{\partial^2 \ell({\bm
\eta}_0)}{\partial{{\bm \eta}}\partial{\bm \eta}'} -
E\frac{\partial^2 \ell({\bm \eta}_0)}{\partial{{\bm \eta}}\partial{\bm
\eta}'}\} = o_p(1)$.  Here  a matrix is said to be $ {\cal O}_p(1)$
(or $o_p(1)$) if all of its elements are of order ${\cal O}_p(1)$
(or $o_p(1)$). The second-order derivatives of $\ell({\bm \eta})$ are
\begin{eqnarray*}
\frac{\partial^2\ell({\bm \eta})}{\partial{{\bm \beta}}\partial{\bm
\beta}'} &=& -(\sigma^2)^{-1}{\bf X}_{nm}'{\bf S}_{nm}'({\bm
\xi}){\bf
S}_{nm}({\bm \xi}){\bf X}_{nm}\\
\frac{\partial^2\ell({\bm \eta})}{\partial{{\bm
\beta}}\partial{\sigma^2}} &=& -(\sigma^4)^{-1}{\bf X}_{nm}'{\bf
S}_{nm}'({\bm \xi}){\bm \nu}_{nm}\\
\frac{\partial^2\ell({\bm \eta})}{\partial{{\bm
\beta}}\partial{\theta_k}} &=& -(\sigma^2)^{-1}{\bf X}_{nm}'\{{\rm
diag}({\bf W}'_{nk}{\bf S}_n({\bm \theta}) + {\bf S}'_n({\bm
\theta}){\bf W}_{nk},\dots,{\bf W}'_{nk}{\bf S}_n({\bm \theta}) +
{\bf S}'_n({\bm \theta}){\bf W}_{nk}) \}{\bf S}^{-1}_{nm}({\bm
\xi}){\bm \nu}_{nm}, k=1,\dots,q\\
\frac{\partial^2\ell({\bm \eta})}{\partial{{\bm
\beta}}\partial{\alpha}} &=& -(\sigma^2)^{-1}{\bf X}_{nm}'\{{\bf
B}_{nm}{\bf B}'_{nm}\}{\bf S}^{-1}_{nm}({\bm \xi}){\bm \nu}_{nm},
k=1,\dots,q
\end{eqnarray*}
where$${\bf B}_{nm} =
\begin{pmatrix}
-\alpha{\bf I}_n & {\bf S}_n({\bm \theta}) &   &   &   &\\
0                      & \ddots & \ddots &   &\\
0                      &        & -\alpha{\bf I}_n & {\bf S}_n({\bm \theta})  &   &\\
0                       & \cdots                  & 0     & 0
\end{pmatrix}
$$
\begin{eqnarray*}
\frac{\partial^2\ell({\bm
\eta})}{\partial{{\sigma^2}}\partial{\sigma^2}} &=&
(2\sigma^6)^{-1}(-2{\bm \nu}'_{nm}{\bm \nu}_{nm}+ nm\sigma^2)\\
\frac{\partial^2\ell({\bm
\eta})}{\partial{{\theta}_k}\partial{\sigma^2}} &=&
-(\sigma^4)^{-1}{\bm \nu}'_{nm}{\rm diag}({\bf W}_{nk}{\bf
S}^{-1}_n({\bm \theta}),\dots,{\bf W}_{nk}{\bf S}^{-1}_n({\bm
\theta})){\bm
\nu}_{nm}\\
\frac{\partial^2\ell({\bm
\eta})}{\partial{{\theta_k}}\partial{\theta_l}} &=& -m {\rm
tr}\{{\bf W}_{nk}{\bf S}^{-1}_n({\bm \theta}){\bf W}_{nl}{\bf
S}^{-1}_n({\bm \theta})\} - (\sigma^2)^{-1}{\bm \nu}'_{nm}{{\bf
S}'}^{-1}_{nm}({\bm \xi}){\rm diag}({\bf W}'_{nl}{\bf
W}_{nk},\dots,{\bf W}'_{nl}{\bf
W}_{nk}){\bf S}^{-1}_{nm}({\bm \xi}){\bm \nu}_{nm}\\
\frac{\partial^2\ell({\bm
\eta})}{\partial{{\alpha}}\partial{\sigma^2}} &=&
-(\sigma^4)^{-1}{\bm \nu}'_{nm}{\bf F}_{nm}{\bf S}^{-1}_{nm}({\bm
\xi}){\bm \nu}_{nm}\\
\frac{\partial^2\ell({\bm
\eta})}{\partial{{\alpha}}\partial{\alpha}} &=& -(\sigma^2)^{-1}{\bm
\nu}'_{nm}{\bf S}^{'-1}_{nm}({\bm \xi}){\rm diag}({\bf
I}_n,\dots,{\bf
I}_n,0){\bf S}^{-1}_{nm}({\bm \xi}){\bm \nu}_{nm}\\
\frac{\partial^2\ell({\bm
\eta})}{\partial{{\alpha}}\partial{\theta_k}} &=& 0
\end{eqnarray*}
By $(A.1)-(A.5)$, we have
\begin{eqnarray*}
&&(nm)^{-1}\{\frac{\partial^2\ell({\tilde{\bm \eta}_{nm}})}{\partial{{\bm \beta}}\partial{\bm \beta}'} - \frac{\partial^2\ell({\bm \eta}_0)}{\partial{{\bm \beta}}\partial{\bm \beta}'}\}\\
&=&(nm)^{-1}\{-(\tilde{\sigma}^2_{nm})^{-1}{\bf X}'_{nm}{\bf S}'_{nm}({\bm \xi}){\bf S}_{nm}({\bm \xi}){\bf X}_{nm}+ ({\sigma^2_0})^{-1}{\bf X}'_{nm}{\bf S}'_{0nm}{\bf S}_{0nm}{\bf X}_{nm}\}\\
&=&(nm)^{-1}{\bf X}'_{nt}{\bf S}'_{0nm}{\bf S}_{0nm}{\bf
X}_{nm}(\frac{1}{{\sigma^2_0}}-\frac{1}{\tilde{\sigma}^2_{nm}}) +
(nm{\tilde{\sigma}^2_{nm}})^{-1}\{{\bf X}'_{nm}{\bf S}'_{0nm}{\bf
S}_{0nm}{\bf X}_{nm} - {\bf X}'_{nm}{\bf S}'_{nm}({\tilde{\bm
\xi}_{nm}}){\bf S}_{nm}({\tilde{\bm \xi}_{nm}}){\bf X}_{nm}\}
\\&=& o_p(1)\\
\\
&&(nm)^{-1}\{\frac{\partial^2\ell({\tilde{\bm \eta}_{nm}})}{\partial{{\bm \beta}}\partial{{\sigma}^2}} - \frac{\partial^2\ell({\bm \eta}_0)}{\partial{{\bm \beta}}\partial{{\sigma}^2}}\}\\
&=&(nm)^{-1}{\bf X}'_{nm}{\bf S}'_{0nm}{\bm \nu}_{0nm}(\frac{1}{{\sigma_0}^4}-\frac{1}{\tilde{\sigma}^4_{nm}})+ (nm{\tilde{\sigma}^4_{nm}})^{-1}{\bf X}'_{nm}{\bf E}_{nm} {\bf S}^{-1}_{0nm}{\bm \nu}_{0nm}\\
&\quad&-({\tilde{\sigma}^4_{nm}})^{-1}{\bf X}'_{nm}{\bf S}'_{nm}({\tilde{\xi}_{nm}}){\bf S}_{nm}({\tilde{\xi}_{nm}}){\bf X}_{nm}({\bm \beta}_0-{\tilde{\bm \beta}_{nm}})\\
&=&o_p(1)
\end{eqnarray*}
where$$ {\bf E}_{nm}=
\begin{pmatrix}
{\bf S}'_{0n}{\bf S}_{0n}-{\bf S}'_n({\tilde{\bm \theta}}_{nm}){\bf S}_n({\tilde{\bm \theta}}_{nm}) &\alpha{\sum^q_{k=1}({\bm \theta}_{0k}-{\tilde{\bm \theta}_{nk}}){\bf W}_{nk}} & & \\
\alpha{\sum^q_{k=1}({\bm \theta}_{0k}-{\tilde{\bm \theta}_{nk}}){\bf W}_{nk}} &\ddots   &\ddots & \\
\ddots & \ddots &{\alpha}{\sum^q_{k=1}({\bm \theta}_{0k}-{\tilde{\bm \theta}_{nk}}){\bf W}_{nk}} &\\
 &\alpha{\sum^q_{k=1}({\bm \theta}_{0k}-{\tilde{\bm \theta}_{nk}}){\bf W}_{nk}}       & {\bf S}'_{0n}{\bf S}_{0n}-{\bf S}'_n({\tilde{\bm \theta}}_{nm}){\bf S}_n({\tilde{\bm
 \theta}}_{nm})&
\end{pmatrix},\quad
{\bf S}_n({\bm \theta}_0) = {\bf S}_{0n}
$$
\begin{eqnarray*}
&&(nm)^{-1}\{\frac{\partial^2\ell({\tilde{\bm \eta}_{nm}})}{\partial{{\bm \beta}}\partial{{\theta}_k}} - \frac{\partial^2\ell({\bm \eta}_0)}{\partial{{\bm \beta}}\partial{{\theta}_k}}\}\\
&=&(nm)^{-1}\{-({\tilde{\sigma}^2_{nm}})^{-1}{\bf X}'_{nm}\{{{\bf J}_{nm}(\tilde{\bm \xi}_{nm})}+{{\bf J}'_{nm}(\tilde{\bm \xi}_{nm})}\}{\bf S}^{-1}_{nm}({\tilde{\bm \xi}_{nm}}){\bm \nu}_{nm}({\tilde{\bm \xi}_{nm}})+({\sigma_0}^2)^{-1}{\bf X}'_{nm}({{\bf J}_{0nm}}+{{\bf J}'_{0nm}}){\bf S}^{-1}_{0nm}{\bm \nu}_{0nm}\}\\
&=&(nm)^{-1}\{-({\tilde{\sigma}^2_{nm}})^{-1}{\bf X}'_{nm}\{{{\bf
J}_{nm}(\tilde{\bm \xi}_{nm})}+{{\bf J}'_{nm}(\tilde{\bm
\xi}_{nm})}\}({\bf Y}_{nm}-{\bf X}_{nm}{\bm \beta}_0)\\
&\quad&+ ({\sigma_0}^2)^{-1}{\bf X}'_{nm}({{\bf J}_{0nm}}+{{\bf
J}'_{0nm}})({\bf Y}_{nm}-{\bf X}_{nm}{\bm \beta}_0)
- ({\tilde{\sigma}^2_{nm}})^{-1}{\bf X}'_{nm}\{{{\bf J}_{nm}(\tilde{\bm \xi}_{nm})}+{{\bf J}'_{nm}(\tilde{\bm \xi}_{nm})}\}{\bf X}_{nm}({\bm \beta}_0 - {\tilde {\bm \beta}_{nm}})\}\\
&=&o_p(1)
\end{eqnarray*}
where ${\bf J}_{nm}={\rm diag}({\bf W}'_{nk}{\bf S}_n({\bm
\theta}),\cdots,{\bf W}'_{nk}{\bf S}_n({\bm \theta}))$
\begin{eqnarray*}
&&(nm)^{-1}\{\frac{\partial^2\ell({\tilde{\bm \eta}_{nm}})}{\partial{{\bm \beta}}\partial{{\alpha}}} - \frac{\partial^2\ell({\bm \eta}_0)}{\partial{{\bm \beta}}\partial{{\alpha}}}\}\\
&=&(nm)^{-1}\{-({\tilde{\sigma}^2_{nm}})^{-1}{\bf X}'_{nm}\{{{\bf B}_{nm}(\tilde{\bm \xi}_{nm})}+{{\bf B}'_{nm}(\tilde{\bm \xi}_{nm})}\}({\bf Y}_{nm}-{\bf X}_{nm}{\tilde{\bm \beta}_{nm}})\\
&\quad& + ({\sigma_0}^2)^{-1}{\bf X}'_{nm}({{\bf B}_{0nm}}+{{\bf B}'_{0nm}})({\bf Y}_{nm}-{\bf X}_{nm}{\bm \beta}_0)\}\\
&=&(nm)^{-1}\{-({\tilde{\sigma}^2_{nm}})^{-1}{\bf X}'_{nm}\{{{\bf B}_{nm}(\tilde{\bm \xi}_{nm})}+{{\bf B}'_{nm}(\tilde{\bm \xi}_{nm})}\}({\bf Y}_{nm}-{\bf X}_{nm}{\bm \beta}_0)\\
&\quad& + ({\sigma_0}^2)^{-1}{\bf X}'_{nm}({{\bf B}_{0nm}}+{{\bf B}'_{0nm}})({\bf Y}_{nm}-{\bf X}_{nm}{\bm \beta}_0)- ({\tilde{\sigma}^2_{nm}})^{-1}{\bf X}'_{nm}\{{{\bf B}_{nm}(\tilde{\bm \xi}_{nm})}+{{\bf B}'_{nm}(\tilde{\bm \xi}_{nm})}\}{\bf X}_{nm}({\bm \beta}_0 -{\tilde{\bm \beta}_{nm}})\}\\
&=& o_p(1)\\
\\
&&(nm)^{-1}\{\frac{\partial^2\ell({\tilde{\bm \eta}_{nm}})}{\partial{{\theta}_k}\partial{{\theta}_l}} - \frac{\partial^2\ell({\bm \eta}_0)}{\partial{{\theta}_k}\partial{{\theta}_l}}\}\\
&=&(nm)^{-1}[-m{\rm tr}\{{\bf W}_{nk}{\tilde{\bf S}^{-1}_n({\bm \theta})}{\bf W}_{nl}{\tilde{\bf S}^{-1}_n({\bm \theta})}\}+m{\rm tr}\{{\bf W}_{nk}{{\bf S}^{-1}_{0n}}{\bf W}_{nl}{{\bf S}^{-1}_{0n}}\}] \\
&\quad&- (nm)^{-1}{\bm \nu}'_{0nm}{{\bf S}'_{0nm}}^{-1}{\rm diag}({\bf W}'_{nl}{\bf W}_{nk},\cdots,{\bf W}'_{nl}{\bf W}_{nk}){\bf S}^{-1}_{0nm}{\bm \nu}_{0nm}(\frac{1}{{\sigma}^2_0} - \frac{1}{{\tilde{\sigma}^2_{nm}}})\\
&\quad&-(nm)^{-1}\{({\tilde{\sigma}^2_{nm}})^{-1}({\bm \beta}_0 - {\tilde{\bm \beta}_{nm}})'{\bf X}'_{nm}{\rm diag}({\bf W}'_{nl}{\bf W}_{nk},\cdots,{\bf W}'_{nl}{\bf W}_{nk}){\bf X}_{nm}({\bm \beta}_0 - {\tilde{\bm \beta}_{nm}})\\
&\quad&-2({\tilde{\sigma}^2_{nm}})^{-1}({\bm \beta}_0 - {\tilde{\bm \beta}_{nm}})'{\bf X}'_{nm}{\rm diag}({\bf W}'_{nl}{\bf W}_{nk},\cdots,{\bf W}'_{nl}){\bf S}^{-1}_{0nm}{\bm \nu}_{0nm}\}\\
&=&o_p(1)\\
\\
&&(nm)^{-1}\{\frac{\partial^2\ell({\tilde{\bm \eta}_{nm}})}{\partial{{\theta}_k}\partial{{\sigma}^2}} - \frac{\partial^2\ell({\bm \eta}_0)}{\partial{{\theta}_k}\partial{{\sigma}^2}}\}\\
&=& (nm)^{-1}\{(({\sigma}^4_0)^{-1}-({\tilde{\sigma}^4_{nm}})^{-1}){\bm \nu}'_{0nm}{\rm diag}({\bf W}_{nk}{\bf S}^{-1}_{0n},\cdots,{\bf W}_{nk}{\bf S}^{-1}_{0n}){\bm \nu}_{0nm}\}\\
&\quad&-(nm)^{-1}({\tilde{\sigma}^4_{nm}})^{-1}({\bm \beta}'_0-{\tilde{\bm \beta}_{nm}})'{\bf X}'_{nm}{\rm diag}({\tilde{\bf S}'_n}{\bf W}_{nk},\cdots,{\tilde{\bf S}'_n}{\bf W}_{nk}){\bf X}_{nm}({\bm \beta}'_0-{\tilde{\bm \beta}_{nm}})\\
&\quad&-(nm)^{-1}({\tilde{\sigma}^4_{nm}})^{-1}({\bm \beta}'_0-{\tilde{\bm \beta}_{nm}})'{\bf X}'_{nm}{\rm diag}({\tilde{\bf S}'_n}{\bf W}_{nk}{\bf S}^{-1}_{0n},\cdots,{\tilde{\bf S}'_{0n}}{\bf W}_{nk}{\bf S}^{-1}_{0n}){\bm \nu}_{0nm}\\
&\quad&-(nm)^{-1}({\tilde{\sigma}^4_{nm}})^{-1}{\bm \nu}'_{0nm}{\rm diag}({{\bf S}^{'-1}_{0n}}{\tilde{\bf S}'_n}{\bf W}_{nk},\cdots,{{\bf S}^{'-1}_{0n}}{\tilde{\bf S}'_{n}}{\bf W}_{nk}){\bf X}_{nm}({\bm \beta}'_0-{\tilde{\bm \beta}_{nm}})\\
&\quad&-(nm)^{-1}({\tilde{\sigma}^4_{nm}})^{-1}{\bm \nu}'_{0nm}{\rm diag}({{\bf S}^{'-1}_{0n}}({\tilde{\bf S}_n}-{\bf S}_{0n}){\bf W}_{nk}{\bf S}^{-1}_{0n},\cdots,{{\bf S}^{'-1}_{0n}}({\tilde{\bf S}_n}-{\bf S}_{0n}){\bf W}_{nk}{\bf S}^{-1}_{0n}){\bm \nu}_{0nm}\\
&=& o_p(1)\\
\\
&&(nm)^{-1}\{\frac{\partial^2\ell({\tilde{\bm \eta}_{nm}})}{\partial{{\sigma}^2}\partial{{\sigma}^2}} - \frac{\partial^2\ell({\bm \eta}_0)}{\partial{{\sigma}^2}\partial{{\sigma}^2}}\}\\
&=&\{(2{\tilde{\sigma}^4_{nm}})^{-1} - (2{\sigma}^4_0)^{-1}\} + (nm)^{-1}{\bm \nu}'_{0nm}{\bm \nu}_{0nm}(1/{\sigma}^6_0-1/{\tilde{\sigma}^6_{nm}})\\
&\quad&-(nm)^{-1}{\bm \nu}'_{0nm}{{\bf S}^{'-1}_{0nm}}\{{\bf S}'_{0nm}{\bf S}_{0nm}-{\bf S}'_{nm}(\tilde{{\bm \xi}_{nm}}){\bf S}_{nm}(\tilde{{\bm \xi}_{nm}})\}{{\bf S}^{-1}_{0nm}}{\bm \nu}_{0nm}\\
&\quad&+({\bf Y}_{nm}-{\bf X}_{nm}{\bm \beta}_0)'{\bf S}'_{nm}(\tilde{{\bm \xi}_{nm}}){\bf S}_{nm}(\tilde{{\bm \xi}_{nm}}){\bf X}_{nm}({\tilde{\bm \beta}_{nm}} - {\bm \beta}_0)\\
&\quad&+({\tilde{\bm \beta}_{nm}} - {\bm \beta}_0)'{\bf X}'_{nm}{\bf S}'_{nm}(\tilde{{\bm \xi}_{nm}}){\bf S}_{nm}(\tilde{{\bm \xi}_{nm}})({\bf Y}_{nm}-{\bf X}_{nm}{\bm \beta}_0)\\
&\quad&-({\tilde{\bm \beta}_{nm}} - {\bm \beta}_0)'{\bf X}'_{nm}{\bf S}'_{nm}(\tilde{{\bm \xi}_{nm}}){\bf S}_{nm}(\tilde{{\bm \xi}_{nm}}){\bf X}_{nm}({\tilde{\bm \beta}_{nm}} - {\bm \beta}_0)\\
&=&o_p(1)\\
\end{eqnarray*}
\begin{eqnarray*}
&&(nm)^{-1}\{\frac{\partial^2\ell({\tilde{\bm \eta}_{nm}})}{\partial{{\alpha}}\partial{{\sigma}^2}} - \frac{\partial^2\ell({\bm \eta}_0)}{\partial{{\alpha}}\partial{{\sigma}^2}}\}\\
&=&(nm)^{-1}\{-({\tilde{\sigma}^4_{nm}})^{-1}{\bm \nu}'_{nm}({\tilde{\bm \xi}_{nm}}){\bf F}_{nm}{\bf S}^{-1}_{nm}({\tilde{\bm \xi}_{nm}}){\bm \nu}_{nm} + ({\sigma}^4_0)^{-1}{\bm \nu}'_{0nm}{\bf F}_{nm}{\bf S}^{-1}_{0nm}{\bm \nu}_{0nm}\}\\
&=&(nm)^{-1}\{[({\tilde{\sigma}^4_{nm}})^{-1}-({\sigma}^4_0)^{-1}]({\bf Y}_{nm}-{\bf X}_{nm}{\bm \beta}_0)'{\rm diag}({\alpha}_0{\bf I}_n,\cdots,{\alpha}_0{\bf I}_n,0)({\bf Y}_{nm}-{\bf X}_{nm}{\bm \beta}_0)\}\\
&\quad&+(nm)^{-1}({\tilde{\sigma}^4_{nm}})^{-1}\{({\bf Y}_{nm}-{\bf X}_{nm}{\tilde{\bm \beta}_{nm}})'{\rm diag}({\tilde{\alpha}_{nm}}{\bf I}_n,\cdots,{\tilde{\alpha}_{nm}}{\bf I}_n,0){\bf X}_{nm}({\bm \beta}_0 - {\tilde{\bm \beta}_{nm}})\\
&\quad&+({\bm \beta}_0 - {\tilde{\bm \beta}_{nm}}){\bf X}'_{nm}{\rm diag}({\tilde{\alpha}_{nm}}{\bf I}_n,\cdots,{\tilde{\alpha}_{nm}}{\bf I}_n,0)({\bf Y}_{nm}-{\bf X}_{nm}{\bm \beta}_0)\\
&\quad&+({\bf Y}_{nm}-{\bf X}_{nm}{\bm \beta}_0)'{\rm diag}(({\tilde{\alpha}_{nm}-{\alpha}_0}){\bf I}_n,\cdots,({\tilde{\alpha}_{nm}-{\alpha}_0}){\bf I}_n,0)({\bf Y}_{nm}-{\bf X}_{nm}{\bm \beta}_0)\}\\
&=&o_p(1)\\
\\
&&(nm)^{-1}\{\frac{\partial^2\ell({\tilde{\bm \eta}_{nm}})}{\partial{{\alpha}}\partial{{\alpha}}} - \frac{\partial^2\ell({\bm \eta}_0)}{\partial{{\alpha}}\partial{{\alpha}}}\}\\
&=&(nm)^{-1}\{[({\sigma}^2_0)^{-1}-({\tilde{\sigma}^2_{nm}})^{-1}]({\bf Y}_{nm}-{\bf X}_{nm}{\bm \beta}_0)'{\rm diag}({\bf I}_n,\cdots,{\bf I}_n,0)({\bf Y}_{nm}-{\bf X}_{nm}{\bm \beta}_0)\}\\
&\quad&+(nm)^{-1}({\tilde{\sigma}^2_{nm}})^{-1}\{({\bf Y}_{nm}-{\bf X}_{nm}{\bm \beta}_0)'{\rm diag}({\bf I}_n,\cdots,{\bf I}_n,0){\bf X}_{nm}({\tilde{\bm \beta}_{nm}}-{\bm \beta}_0)\\
&\quad&+({\tilde{\bm \beta}_{nm}}-{\bm \beta}_0)'{\bf X}'_{nm}{\rm diag}({\bf I}_n,\cdots,{\bf I}_n,0)({\bf Y}_{nm}-{\bf X}_{nm}{\bm \beta}_0)\\
&\quad&-({\tilde{\bm \beta}_{nm}}-{\bm \beta}_0)'{\bf X}'_{nm}{\rm diag}({\bf I}_n,\cdots,{\bf I}_n,0){\bf X}_{nm}({\tilde{\bm \beta}_{nm}}-{\bm \beta}_0)\}\\
&=&o_p(1)\\
\\
&&(nm)^{-1}\{\frac{\partial^2\ell({\tilde{\bm
\eta}_{nm}})}{\partial{{\alpha}}\partial{{\bm \theta}_k}} -
\frac{\partial^2\ell({\bm \eta}_0)}{\partial{{\alpha}}\partial{{\bm
\theta}_k}}\} = o_p(1)
\end{eqnarray*}
\noindent Thus, $$(nm)^{-1}\{\frac{\partial^2\ell({\tilde{\bm
\eta}_{nm}})}{\partial{{\bm \eta}}\partial{{\bm \eta}'}} -
\frac{\partial^2\ell({\bm \eta}_0)}{\partial{{\bm
\eta}}\partial{{\bm \eta}'}}\}=o_p(1)$$ Further,
under$(A.1)-(A.5)$,we have
\begin{eqnarray*}
(nm)^{-1}\{\frac{\partial^2\ell({\bm \eta}_0)}{\partial{{\bm
\beta}}\partial{{\bm \beta}'}} - E\frac{\partial^2\ell({\bm
\eta}_0)}{\partial{{\bm \beta}}\partial{{\bm \beta}'}}\}&=& 0\\
(nm)^{-1}\{\frac{\partial^2\ell({\bm \eta}_0)}{\partial{{\bm
\beta}}\partial{{\sigma}^2}} - E\frac{\partial^2\ell({\bm
\eta}_0)}{\partial{{\bm \beta}}\partial{{\sigma}^2}}\}&=&
-(nm)^{-1}({\sigma}^4_0)^{-1}{\bf X}'_{nm}{\bf S}'_{0nm}{\bm
\nu}'_{0nm} = (nm)^{-1/2} \times O_p(1) = o_p(1)\\
(nm)^{-1}\{\frac{\partial^2\ell({\bm \eta}_0)}{\partial{{\bm
\beta}}\partial{{\bm \theta}_k}} - E\frac{\partial^2\ell({\bm
\eta}_0)}{\partial{{\bm \beta}}\partial{{\bm \theta}_k}}\} &=&
-(nm)^{-1}({\sigma}^2_0)^{-1}{\bf X}'_{nm}({\bf J}_{0nm}+{\bf
J}'_{0nm}){\bf S}^{-1}_{0nm}{\bm \nu}_{0nm} =  o_p(1)\\
(nm)^{-1}\{\frac{\partial^2\ell({\bm \eta}_0)}{\partial{{\bm
\beta}}\partial{{\alpha}}} - E\frac{\partial^2\ell({\bm
\eta}_0)}{\partial{{\bm \beta}}\partial{{\alpha}}}\}&=&
-(nm)^{-1}({\sigma}^2_0)^{-1}{\bf X}'_{nm}({\bf B}_{0nm}+{\bf
B}'_{0nm}){\bf S}^{-1}_{0nm}{\bm \nu}_{0nm} = o_p(1)\\
(nm)^{-1}\{\frac{\partial^2\ell({\bm
\eta}_0)}{\partial{{\sigma}^2}\partial{{\sigma}^2}} -
E\frac{\partial^2\ell({\bm
\eta}_0)}{\partial{{\sigma}^2}\partial{{\sigma}^2}}\}&=&-(nm)^{-1}({\sigma}^6_0)^{-1}({\bm
\nu}'_{0nm}{\bm \nu}_{0nm} - nm)=o_p(1)\\
(nm)^{-1}\{\frac{\partial^2\ell({\bm \eta}_0)}{\partial{{\bm
\theta}_k}\partial{{\sigma}^2}} - E\frac{\partial^2\ell({\bm
\eta}_0)}{\partial{{\bm \theta}_k}\partial{{\sigma}^2}}\}
 &=& -(nm)^{-1}({\sigma}^4_0)^{-1}\{{\bm \nu}'_{0nm}{\rm diag}({\bf W}_{nk}{\bf S}^{-1}_{0n},\cdots,{\bf W}_{nk}{\bf S}^{-1}_{0n}){\bm \nu}_{0nm} \\ &\quad &- ({\sigma_0}^2)\times m{\rm tr}({\bf W}_{nk}{\bf S}^{-1}_{0n})\} =
 o_p(1)
 \end{eqnarray*}
\begin{eqnarray*}
 (nm)^{-1}\{\frac{\partial^2\ell({\bm \eta}_0)}{\partial{{\bm \theta}_k}\partial{{\bm \theta}_l}} - E\frac{\partial^2\ell({\bm \eta}_0)}{\partial{{\bm \theta}_k}\partial{{\bm \theta}_l}}\}
&=&-(nm)^{-1}({\sigma}^2_0)^{-1}\{{\bm \nu}'_{0nm}{{\bf S}^{'-1}_{0nm}}{\rm diag}({\bf W}'_{nl}{\bf W}_{nk},\cdots,{\bf W}'_{nl}{\bf W}_{nk}){\bf S}^{-1}_{0nm}{\bm \nu}_{0nm}\\
&\quad&- {\sigma}^2_0 \times {\rm tr} \{{{\bf S}^{'-1}_{0nm}}{\rm
diag}({\bf W}'_{nl}{\bf W}_{nk},\cdots,{\bf W}'_{nl}{\bf
W}_{nk}){\bf S}^{-1}_{0nm}\}\} = o_p(1)\\
(nm)^{-1}\{\frac{\partial^2\ell({\bm
\eta}_0)}{\partial{{\alpha}}\partial{{\sigma}^2}} -
E\frac{\partial^2\ell({\bm
\eta}_0)}{\partial{{\alpha}}\partial{{\sigma}^2}}\}&=&
-(nm)^{-1}({\sigma}^4_0)^{-1}\{{\bm \nu}'_{0nm}{\bf F}_{nm}{\bf
S}^{-1}_{0nm}{\bm \nu}_{0nm} - {\sigma}^2_0 \times {\rm tr}\{{\bf
F}_{nm}{\bf S}^{-1}_{0nm}\}\} = o_p(1)\\
(nm)^{-1}\{\frac{\partial^2\ell({\bm
\eta}_0)}{\partial{{\alpha}}\partial{{\alpha}}} -
E\frac{\partial^2\ell({\bm
\eta}_0)}{\partial{{\alpha}}\partial{{\alpha}}}\}&=&
-(nm)^{-1}({\sigma}^2_0)^{-1}\{{\bm \nu}'_{0nm}{{\bf
S}^{'-1}_{0nm}}{\rm diag}({\bf I}_n,\cdots,{\bf I}_n,0){\bf
S}^{-1}_{0nm}{\bm \nu}_{0nm}\\&\quad&- {\sigma}^2_0 \times {\rm
tr}\{{{\bf S}^{'-1}_{0nm}}{\rm diag}({\bf I}_n,\cdots,{\bf
I}_n,0){\bf S}^{-1}_{0nm}\}\}= o_p(1)\\
(nm)^{-1}\{\frac{\partial^2\ell({\bm
\eta}_0)}{\partial{{\alpha}}\partial{{\bm \theta}_k}} -
E\frac{\partial^2\ell({\bm \eta}_0)}{\partial{{\alpha}}\partial{{\bm
\theta}_k}}\} &=& 0
\end{eqnarray*}
Furthermore, the first-order derivatives of
$\ell({\bm \eta})$ at ${\bm \eta}_0$ are linear or
quadratic forms of ${\bm \nu}_{0nm}$ since
\[ \frac{\partial \ell({\bm \eta}_0)}{\partial {\bm \beta}} = ({\sigma}^2_0)^{-1}{\bm X}'_{nm}{\bm S}'_{0nm}{\bm \nu}_{0nm}, \quad \frac{\partial \ell({\bm \eta}_0)}{\partial {\sigma}^2} = (2{\sigma}^4_0)^{-1}({\bm \nu}'_{0nm}{\bm \nu}_{0nm}- nm {\sigma}^2_0)\]
\[ \frac{\partial \ell({\bm \eta}_0)}{\partial \theta_k} = -{\rm tr}({\bm G}_k) + ({\sigma}^2_0)^{-1}{\bm \nu}'_{0nm}{\bm G}_k{\bm \nu}_{0nm}, \quad \frac{\partial \ell({\bm \eta}_0)}{\partial {\alpha}} = ({\sigma}^2_0)^{-1}{\bm \nu}_{0nm}{\bm F}_{nm}{\bm S}^{-1}_{0nm}{\bm \nu}_{0nm} \]
where ${\bm G}_k = {\rm diag}(\bm{I}_m \otimes {\bm W}_{nk}{\bm S}^{-1}_{n}(\bm{\theta}_0))$ ,for $k=1,\ldots,q$.

By (A.5),
\[ (nm)^{-1/2}\frac{\partial \ell({\bm \eta}_0)}{\partial {\bm \beta}} \stackrel{\mathcal{D}}{\rightarrow} N(0,\lim_{nm \rightarrow \infty}(nm)^{-1}({\sigma}^2_0)^{-1}{\bm X}'_{nm}{\bm S}'_{0nm}{\bm S}_{0nm}{\bm X}_{nm}). \]
By a classic central limit theorem,
\[(nm)^{-1/2}\frac{\partial \ell({\bm \eta}_0)}{\partial {\sigma}^2} \stackrel{\mathcal{D}}{\rightarrow} N(0,(2{\sigma}^4_0)^{-1})\]
For asymptotic normality of $(nm)^{-1/2}\frac{\partial \ell({\bm \eta}_0)}{\partial \theta_k}$, (A.2) and (A.4) ensure that ${\bm G}_k$ is uniformly bounded in matrix norm ${||\cdot||}_1$ and ${||\cdot||}_{\infty}$ and the positive definiteness of $\bm{\Sigma}_{\bm{\eta}_0}^{-1}$ ensures that $(nm)^{-1} {\rm Var}(\frac{\partial \ell({\bm \eta}_0)}{\partial \theta_k}) = (nm)^{-1}{\rm tr}({\bm G}^2_k + {\bm G}'_k{\bm G}_k)$ is bounded away from 0. By a central limit theorem for linear-quadratic forms (Theorem 1, \cite*{kelejianp01}), we have
$$(nm)^{-1/2}\frac{\partial \ell({\bm \eta}_0)}{\partial \theta_k} \stackrel{\mathcal{D}}{\rightarrow} N(0,\lim_{nm \rightarrow \infty}(nm)^{-1}{\rm tr}({\bm G}^2_k + {\bm G}'_k{\bm G}_k))$$ for $k = 1,\cdots,q$. Similarly, we can get the asymptotic normality of $(nm)^{-1/2}\frac{\partial \ell({\bm \eta}_0)}{\partial {\alpha}}$
\[ (nm)^{-1/2}\frac{\partial \ell({\bm \eta}_0)}{\partial {\alpha}} \stackrel{\mathcal{D}}{\rightarrow} N(0,\lim_{nm \rightarrow \infty}(nm)^{-1}{\rm tr}({\bm H}^2_{nm} + {\bm H}'_{nm}{\bm H}_{nm}))\]
By Cramer-Wold Theorem and the fact that $\frac{\partial \ell({\bm \eta}_0)}{\partial {\bm \beta}}$ is asymptotically independent of $\frac{\partial \ell({\bm \eta}_0)}{\partial {\alpha}},\frac{\partial \ell({\bm \eta}_0)}{\partial {\sigma}^2}$ and $\frac{\partial \ell({\bm \eta}_0)}{\partial \theta_k}$, for $\hspace{0.1cm} k=1,\cdots,q$, we have
$$(nm)^{-1/2}\frac{\partial \ell({\bm \eta}_0)}{\partial {\bm \eta}} \stackrel{\mathcal{D}}{\rightarrow} N(0,{\Sigma}^{-1}_{{\bm \eta}_0})$$
where ${\Sigma}^{-1}_{{\bm \eta}_0} = \lim_{nm \rightarrow \infty} E\left(-(nm)^{-1}\frac{\partial^2 \ell({\bm \eta}_0)}{\partial{{\bm \eta}}\partial{{\bm \eta}'}}\right)$ and
\begin{eqnarray*}
E\left(-(nm)^{-1}\frac{\partial^2 \ell({\bm \eta}_0)}{\partial{{\bm \beta}}\partial{{\bm \beta}'}}\right)&=& (nm{\sigma}^2_0)^{-1}{\bm X}'_{nm}{\bm S}'_{0nm}{\bm S}_{0nm}{\bm X}_{nm}, \quad E\left(-(nm)^{-1}\frac{\partial^2 \ell({\bm \eta}_0)}{\partial{{\bm \beta}}\partial{{\bm \xi}'}}\right)=0,\\
E\left(-(nm)^{-1}\frac{\partial^2 \ell({\bm \eta}_0)}{\partial{{\bm \beta}}\partial{\sigma^2}}\right)&=&0,\quad E\left(-(nm)^{-1}\frac{\partial^2 \ell({\bm \eta}_0)}{\partial{{\bm \theta}}\partial{\alpha}}\right)=0\\
E\left(-(nm)^{-1}\frac{\partial^2 \ell({\bm \eta}_0)}{\partial{\theta_k}\partial{\theta_l}}\right)&=& (nm)^{-1}{\rm tr}({\bm G}_{k}{\bm G}_{l}+ {\bm G}'_{l}{\bm G}_{k}), k,l=1,\ldots,q\\
E\left(-(nm)^{-1}\frac{\partial^2 \ell({\bm \eta}_0)}{\partial{\theta_k}\partial{\sigma^2}}\right)&=&(nm{\sigma}^2_0)^{-1}{\rm tr}({\bm G}_{k}), k=1,\ldots,q
\end{eqnarray*}
\begin{eqnarray*}
E\left(-(nm)^{-1}\frac{\partial^2 \ell({\bm \eta}_0)}{\partial{\alpha^2}}\right)&=&(nm)^{-1}{\rm tr}({\bm H}^2_{nm} + {\bm H}'_{nm}{\bm H}_{nm}), \\
E\left(-(nm)^{-1}\frac{\partial^2 \ell({\bm \eta}_0)}{\partial{\alpha}\partial{\sigma^2}}\right)&=&(nm{\sigma}^2_0)^{-1}{\rm tr}({\bm H}_{nm}),\quad
E\left(-(nm)^{-1}\frac{\partial^2 \ell({\bm \eta}_0)}{\partial{(\sigma^2)^2}}\right)=(2{\sigma}^4_0)^{-1}
\end{eqnarray*}
Assumptions (A.1)-(A.5) ensure that all the elements in ${\Sigma}^{-1}_{{\bm \eta}_0}$ exist.  Then it follows that
$$(nm)^{1/2}(\hat{\bm{\eta}}_{nm} - {\bm \eta}_0) = -\left\{(nm)^{-1}\frac{\partial^2 \ell({\bm \eta}_0)}{\partial{{\bm \eta}}\partial{{\bm \eta}'}} + o_p(1)\right\}^{-1}(nm)^{-1/2}\frac{\partial \ell({\bm \eta}_0)}{\partial {\bm \eta}} \stackrel{\mathcal{D}}{\rightarrow} N(0,{\Sigma}_{{\bm \eta}_0}).$$
and thus the result of this theorem holds. 

In this theorem, we assume the existence and positive definiteness of the covariance matrix $\bm{\Sigma}_{\bm{\eta}_0}$. Here we discuss a simple example about the validation of the existence of $\bm{\Sigma}_{\bm{\eta}_0}$. Suppose we only consider the first-order spatial neighborhood, i.e. $q=1$. Let $\bm{\tau} = (\bm{\theta}', \sigma^2)'$. Then 
\begin{eqnarray*}
&&{\bm \Sigma}^{-1}_{\eta_0} = \lim_{nm \rightarrow \infty} E\left(-(nm)^{-1}\frac{\partial^2 \ell({\bm \eta}_0)}{\partial{{\bm \eta}}\partial{{\bm \eta}'}}\right)\\ 
&=& \lim_{nm \rightarrow \infty} \left(\begin{array}{cccc}
(nm\sigma^2_0)^{-1}{\bm X}'_{nm}{\bm S}'_{0nm}{\bm S}_{0nm}{\bm
X}_{nm} & 0  & 0\\
0 & E\left(-(nm)^{-1}\frac{\partial^2 \ell({\bm
\eta}_0)}{\partial{{\bm \tau}}\partial{{\bm \tau}'}}\right) & 0\\
0 & 0 & (nm)^{-1}{\rm tr}({\bm H}^2_{nm} + {\bm H}'_{nm}{\bm
H}_{nm})
\end{array}\right)\\
&=& \lim_{nm \rightarrow \infty} (nm)^{-1}
\left(\begin{array}{cccc}
(\sigma^2_0)^{-1}{\bm X}'_{nm}{\bm S}'_{0nm}{\bm S}_{0nm}{\bm
X}_{nm} & 0 & 0 & 0\\
0 & {\rm tr}({\bm G}^2_{1}+ {\bm G}'_{1}{\bm G}_{1}) & (\sigma^2_0)^{-1}{\rm tr}({\bm G}_{1}) & 0\\
0 & (\sigma^2_0)^{-1}{\rm tr}({\bm G}_{1}) & nm(2\sigma^4_0)^{-1} &
0 \\
0 & 0 & 0 & {\rm tr}({\bm H}^2_{nm} + {\bm H}'_{nm}{\bm
H}_{nm}),
\end{array}\right)
\end{eqnarray*}
which is a block diagonal matrix. Assumption~(A.5) guarantees the existence and nonsingularity of $\lim_{nm \rightarrow \infty} (nm\sigma^2_0)^{-1}{\bm X}'_{nm}{\bm S}'_{0nm}{\bm S}_{0nm}{\bm X}_{nm}$.  Assumptions (A.2)--(A.4) ensure that $ \lim_{nm \rightarrow \infty} (nm)^{-1}{\rm tr}({\bm H}^2_{nm} + {\bm H}'_{nm}{\bm
H}_{nm})$ is bounded and bounded away from zero.  For the middle block, we have
\begin{eqnarray*} 
\left| E\left(-(nm)^{-1}\frac{\partial^2 \ell({\bm
\eta}_0)}{\partial{{\bm \tau}}\partial{{\bm \tau}'}}\right) \right| &=&  (2nm \sigma^4_0)^{-1}\left[ {\rm tr}({\bm G}^2_{1}+ {\bm
G}'_{1}{\bm G}_{1}) - 2(nm)^{-1}{\rm tr}^2({\bm G}_1) \right] \\
&=& (4nm\sigma^4_0)^{-1}{\rm tr}(\bm{G}^{b}\bm{G}^b)=(4nm\sigma^4_0)^{-1} ||\bm{G}^b||^2_{\rm F}  > 0,
\end{eqnarray*}
where $\bm{G}^b = \bm{G}_1+\bm{G}_1'-2(nm)^{-1}{\rm tr}(\bm{G}_1)\bm{I}_{nm}$ and $||\cdot||_{\rm F}$ denote the Frobenius norm.  Hence, the matrix $E\left(-(nm)^{-1}\frac{\partial^2 \ell({\bm
\eta}_0)}{\partial{{\bm \eta}}\partial{{\bm \eta}'}}\right)$ is positive definite in the large-sample case. Furthermore, assumption~(A.6) guarantees that it is positive definite in the limit (see, e.g. \cite*{lee04}).

%%%%%%%%%%%%%%%%%%%%%%
% h_n go to infinity %
%%%%%%%%%%%%%%%%%%%%%%

\subsection*{The case $h_n \rightarrow \infty$}

Now, we establish the asymptotic normality of the MLE $\hat{\bm{\theta}}_{nm}$ by showing the asymptotic normality of $\left\{h_n/(nm)\right\}^{1/2}\frac{\partial \ell(\bm{\xi}_0)}{\partial \bm{\theta}}$ and the convergence in probability of $h_n (nm)^{-1}\frac{\partial^2 \ell(\tilde{\bm{\xi}}_{nm})}{\partial \bm{\theta} \partial \bm{\theta}'}$, where $\tilde{\bm{\xi}}_{nm}=\lambda\bm{\xi}_0+(1-\lambda)\hat{\bm{\xi}}_{nm}$ for $\lambda\in(0,1)$ converges to $\bm{\xi}_0$ in probability. The asymptotic normality of the MLE $\hat{\alpha}_{nm}$ will be established by showing the asymptotic normality of $(nm)^{-1/2}\frac{\partial \ell(\bm{\xi}_0)}{\partial \alpha}$ and the convergence in probability of $(nm)^{-1}\frac{\partial^2 \ell(\tilde{\bm{\xi}}_{nm})}{\partial \alpha^2}$.

For convergence of $h_n (nm)^{-1}\frac{\partial^2 \ell(\tilde{\bm{\xi}}_{nm})}{\partial \bm{\theta} \partial \bm{\theta}'}$, we show that $h_n (nm)^{-1}\{\frac{\partial^2 \ell(\tilde{\bm{\xi}}_{nm})}{\partial\bm{\theta}\partial\bm{\theta}'} - \frac{\partial^2 \ell(\bm{\xi}_0)}{\partial\bm{\theta}\partial\bm{\theta}'}\}=o_p(1)$ and $h_n (nm)^{-1}\left\{ \frac{\partial^2 \ell(\bm{\xi}_0)}{\partial \bm{\theta} \partial \bm{\theta}'} - E \frac{\partial^2 \ell(\bm{\xi}_0)}{\partial \bm{\theta} \partial \bm{\theta}'} \right\} = o_p(1)$, under (A.1)--(A.5). On the other hand, we have that $(nm)^{-1}\{\frac{\partial^2 \ell(\tilde{\bm{\xi}}_{nm})}{\partial\alpha^2} - \frac{\partial^2 \ell(\bm{\xi}_0)}{\partial\alpha^2}\}=o_p(1)$ and $(nm)^{-1}\left\{ \frac{\partial^2 \ell(\bm{\xi}_0)}{\partial \alpha^2} - E \frac{\partial^2 \ell(\bm{\xi}_0)}{\partial \alpha^2} \right\} = o_p(1)$.  Under(A.1)-(A.5), We note that
\begin{eqnarray*}
\frac{\partial{\ell({\bm \xi})}}{\partial{\theta_1}} &=  &
-\{2{\hat{\sigma}}^2_{nm}({\bm \xi})\}^{-1}{\bm
\nu}'_{0nm}\frac{\partial{\bf B}_{nm}({\bm
\xi})}{\partial{\theta_1}}{\bm \nu}_{0nm} - {\rm tr}\{({\bf
I}_m\otimes{\bf W}_{n1}){\bf S}^{-1}_{nm}({\bm \xi} )\}\\
& = &\{{\hat{\sigma}}^2_{nm}({\bm \xi})\}^{-1} {\bm \nu}'_{0nm}{\bf
T}_{nm,1}({\bm \xi}){\bm \nu}_{0nm} - {\rm tr} \{({\bf I}_m
\otimes{\bf W}_{n1}){\bf S}^{-1}_{nm}({\bm \xi} )\}\\
\frac{\partial{\ell^2({\bm \xi})}}{\partial{{\theta_1}^2}}
&=&2\{nm\hat{\sigma}^4_{nm}(\bm \xi)\}^{-1}\{{\bm \nu}'_{0nm}{\bf
T}_{nm,1}({\bm \xi}){\bm \nu}_{0nm}\}^2 + \{\hat{\sigma}^2_{nm}(\bm
\xi)\}^{-1}{\bm \nu}'_{0nm}\frac{\partial{\bf T}_{nm,1}({\bm
\xi})}{\partial{\theta_1}}{\bm \nu}_{0nm} \\&\quad&-{\rm tr}\{{\bf
S}^{-1}_{nm}({\bm \xi})({\bf I}_m\otimes{\bf W}_{n1}){\bf
S}^{-1}_{nm}({\bm \xi})({\bf I}_m\otimes{\bf W}_{n1})\}
\end{eqnarray*}
\begin{eqnarray*}
&&\bm{T}_{nm,1}(\bm{\xi})=\\
&&{{\bm S}^{'-1}_{nm}}({\bm \xi})\left[({\bm I}_m \otimes{\bm W}_{n1}) {\bm M}_{nm}({\bm \xi}){\bm S}_{nm}({\bm \xi}) - {\bm S}'_{nm}({\bm \xi})({\bm I}_m\otimes{\bm W}_{n1}){\bm X}_{nm}\{{\bm X}'_{nm}{\bm S}'_{nm}({\bm \xi}){\bm S}_{nm}({\bm \xi}){\bm X}_{nm}\}^{-1}{\bm X}'_{nm}{\bm S}'_{nm}({\bm \xi}){\bm S}_{nm}({\bm \xi})\right.\\
&+& {\bm S}'_{nm}({\bm \xi}){\bm S}_{nm}({\bm \xi}){\bm X}_{nm}\{{\bm X}'_{nm}{\bm S}'_{nm}({\bm \xi}){\bm S}_{nm}({\bm \xi}){\bm X}_{nm}\}^{-1}{\bm X}'_{nm} ({\bm I}_m\otimes{\bm W}_{n1}) {\bm S}_{nm}({\bm \xi})({\bm \xi}){\bm X}_{nm} \{{\bm X}'_{nm}{\bm S}'_{nm}({\bm \xi}){\bm S}_{nm}({\bm \xi}){\bm X}_{nm}\}^{-1}\\
&&\quad\quad\quad\quad\quad\quad\quad\quad\quad\left.{\bm
X}'_{nm}{\bm S}'_{nm}({\bm \xi}){\bm S}_{nm}({\bm \xi})\right]{\bm
S}^{-1}_{nm}({\bm \xi})
\end{eqnarray*}
 with ${\bm M}_{nm}({\bm \xi}) = {\bm I}_{nm} - {\bm S}_{nm}({\bm
\xi}){\bm X}_{nm}\{{\bm X}'_{nm}{\bm S}'_{nm}({\bm \xi}){\bm
S}_{nm}({\bm \xi}){\bm X}_{nm}\}^{-1}{\bm X}'_{nm}{\bm S}'_{nm}({\bm
\xi})$. By(A.1)-(A.5)
\begin{eqnarray*}
h_n(nm)^{-1}{\bm \nu}'_{0nm}{\bf T}_{nm,1}({\bm \xi}){\bm \nu}_{0nm}
&=& O_p(1)\\
h_n(nm)^{-1}{\bm \nu}'_{0nm}\frac{\partial{\bf T}_{nm,1}({\bm
\xi})}{\partial{\theta_1}}{\bm \nu}_{0nm} &=&
h_n(nm)^{-1}{\sigma}^2_0 {\rm tr}\{\frac{\partial{\bf T}_{nm,1}({\bm
\xi})}{\partial{\theta_1}}\} + o_p(1)
\end{eqnarray*}

Since ${\bf T}_{nm,1}$ and $\frac{\partial{\bf
T}_{nm,1}({\bm \xi})}{\partial{\theta_1}}$ are uniformly bounded in
either matrix norm $||\cdot||_1$ or $||\cdot||_{\infty}$, thus,
under (A.1) - (A.5),
\begin{eqnarray*}
\frac{\partial{\ell^2({\bm \xi})}}{\partial{{\theta_1}^2}} &=
&2\{nm\hat{\sigma}^4_{nm}(\bm \xi)\}^{-1}\{{\bm \nu}'_{0nm}{\bf
T}_{nm,1}({\bm \xi}){\bm \nu}_{0nm}\}^2 + \{\hat{\sigma}^2_{nm}(\bm
\xi)\}^{-1}{\bm \nu}'_{0nm}\frac{\partial{\bf T}_{nm,1}({\bm
\xi})}{\partial{\theta_1}}{\bm \nu}_{0nm} \\&\quad&-
h_n(nm)^{-1}{\rm tr}\{({\bf I}_m\otimes{\bf W}_{n1}){\bf
S}^{-1}_{nm}({\bm \xi} )({\bf I}_m\otimes{\bf W}_{n1}){\bf
S}^{-1}_{nm}({\bm \xi} )\}
\\&=& h_n(nm)^{-1}{\rm tr}\{\frac{\partial{\bf T}_{nm,1}({\bm
\xi})}{\partial{\theta_1}}\} - h_n(nm)^{-1}{\rm tr}\{({\bf
I}_m\otimes{\bf W}_{n1}){\bf S}^{-1}_{nm}({\bm \xi})({\bf
I}_m\otimes{\bf W}_{n1}){\bf S}^{-1}_{nm}({\bm \xi} )\} + o_p(1)
\end{eqnarray*}
and for $\tilde {\bm \xi}_{nm} = \lambda{\bm \xi}_0 +
(1-\lambda)\hat{{\bm \xi}}_{nm}$,
\begin{align*}
&h_n(nm)^{-1}\{\frac{\partial{\ell^2(\tilde{{\bm
\xi}}_{nm})}}{\partial{{\theta_1}^2}} - \frac{\partial{\ell^2({\bm
\xi}_0)}}{\partial{{\theta_1}^2}}\}
\\ = & h_n(nm)^{-1}\{{\rm tr}\{\frac{\partial{\bf T}_{nm,1}(\tilde{{\bm
\xi}}_{nm})}{\partial{\theta_1}}\} - {\rm tr}\{\frac{\partial{\bf
T}_{nm,1}({\bm \xi}_0)}{\partial{\theta_1}}\}\} -h_n(nm)^{-1}[{\rm
tr}\{({\bf I}_m\otimes{\bf W}_{n1}){\bf S}^{-1}_{nm}(\tilde{{\bm
\xi}}_{nm}) ({\bf I}_m\otimes{\bf W}_{n1}){\bf
S}^{-1}_{nm}(\tilde{{\bm \xi}}_{nm})\}
\\ &- {\rm tr}\{({\bf
I}_m\otimes{\bf W}_{n1}){\bf S}^{-1}_{0nm}({\bf I}_m\otimes{\bf
W}_{n1}){\bf S}^{-1}_{0nm}\}] + o_p(1)
\\= & (\tilde{{\bm
\xi}}_{nm} - {\bm \xi}_0)'\mathcal{O}(1) - (\tilde{{\bm \xi}}_{nm} -
{\bm \xi}_0)'\mathcal{O}(1) +o_p(1)\\
 =& o_p(1)
\end{align*}

By similar argument for $\theta_k$ and $\theta_{k'}$,$k,k'
= 1,\dots,q$, we have
$$h_n(nm)^{-1}\{\frac{\partial^2{\ell(\tilde{{\bm
\xi}}_{nm})}}{\partial{\bm \theta}\partial {\bm \theta'}} -
\frac{\partial^2{\ell({\bm \xi}_0)}}{\partial{\bm \theta}\partial
{\bm \theta'}}\} = o_p(1)$$
Furthermore,
\begin{align*}
-h_n(nm)^{-1}\frac{\partial{\ell^2({\bm
\xi}_0)}}{\partial{{\theta_1}^2}} &= -h_n(nm)^{-1}
tr\{\frac{\partial{\bf T}_{nm,1}({\bm \xi}_0)}{\partial{\theta_1}}\}
+ h_n(nm)^{-1} tr({\bf G}^2_1) +o_p(1)
\\&= -h_n(nm)^{-1}\{tr({\bf G}'_1{\bf G}_1)+ tr({\bf G}^2_1)\}
+o_p(1)
\end{align*}
where ${\bf G}_1 = ({\bf I}_m\otimes{\bf W}_{n1}){\bf
S}^{-1}_{0nm}$. 

Also, we have $h_n(nm)^{-1}\{\frac{\partial^2{\ell({\bm
\xi}_0})}{\partial{\bm \theta}\partial {\bm \theta'}} -E
\frac{\partial^2{\ell({\bm \eta}_0)}}{\partial{\bm \theta}\partial
{\bm \theta'}}\} = o_p(1)$. Thus,
$h_n(nm)^{-1}\frac{\partial^2{\ell(\tilde{{\bm
\xi}}_{nm})}}{\partial{\bm \theta}\partial {\bm \theta'}}
\rightarrow \Sigma^{-1}_{{\bm \theta}_0} = \lim_{nm \rightarrow
\infty}h_n(nm)^{-1} E \frac{\partial^2{\ell({\bm
\eta}_0)}}{\partial{\bm \theta}\partial {\bm \theta'}}$.

To establish the asymptotic
normality of
$\left\{h_n/(nm)\right\}^{1/2}\frac{\partial
\ell(\bm{\xi}_0)}{\partial \bm{\theta}}$, we apply the central limit
theorem for linear-quadratic forms in Appendix A of \cite{lee04}.
Note that
\begin{eqnarray*}
\left\{h_n/(nm)\right\}^{1/2}\frac{\partial{\ell({\bm \xi}_0)}}{\partial{\theta_1}} & =&  \left\{h_n/(nm)\right\}^{1/2} \left[\{\hat{\sigma}^2_{nm}({\bm \theta}_0)\}^{-1}{\bm \nu}'_{0nm}{\bm T}_{nm,1}({\bm \xi}_0){\bm \nu}_{0nm} - {\rm tr}\left\{({\bm I}_m \otimes{\bm W}_{n,1}){\bm S}^{-1}_{0nm}\right\}\right]\\
& = & \{\hat{\sigma}^2_{nm}({\bm
\theta}_0)\}^{-1}\left\{h_n/(nm)\right\}^{1/2}\{{\bm \nu}'_{0nm}{\bm
G}'_1{\bm \nu}_{0nm} - {\sigma}^2_0{\rm tr}({\bm G}_1)\} + o_p(1)
\end{eqnarray*}
where
\begin{eqnarray*}
&&\bm{T}_{nm,1}(\bm{\xi}_0)\\
&=&{{\bm S}^{-1}_{0nm}}'\left[({\bm I}_m \otimes{\bm W}_{n1}) {\bm M}_{0nm}{\bm S}_{0nm} - {\bm S}'_{0nm}({\bm I}_m\otimes{{\bm W}_{n1}}){\bm X}_{nm}\{{\bm X}'_{nm}{\bm S}'_{0nm}{\bm S}_{0nm}{\bm X}_{nm}\}^{-1}{\bm X}'_{nm}{\bm S}'_{0nm}{\bm S}_{0nm}\right.\\
&+& {\bm S}'_{0nm}{\bm S}_{0nm}{\bm X}_{nm}\{{\bm X}'_{nm}{\bm S}'_{0nm}{\bm S}_{0nm}{\bm X}_{nm}\}^{-1}{\bm X}'_{nm} ({\bm I}_m\otimes{\bm W}_{n,1}) {\bm S}_{0nm}{\bm X}_{nm} \{{\bm X}'_{nm}{\bm S}'_{0nm}{\bm S}_{0nm}{\bm X}_{nm}\}^{-1}\\
&&\quad\quad\quad\quad\quad\left.{\bm X}'_{nm}{\bm S}'_{0nm}{\bm
S}_{0nm}\right]{\bm S}^{-1}_{0nm}
\end{eqnarray*}
with ${\bm M}_{0nm} = {\bm I}_{nm} - {\bm S}_{0nm}{\bm X}_{nm}\{{\bm
X}'_{nm}{\bm S}'_{0nm}{\bm S}_{0nm}{\bm X}_{nm}\}^{-1}{\bm
X}'_{nm}{\bm S}'_{0nm}$. (A.2) and (A.4) ensure that $\bm{G}_1$ is
uniformly bounded in matrix norms $||\cdot||_{1}$ and
$||\cdot||_{\infty}$ and the positive definiteness of
$\bm{\Sigma}^{-1}_{\bm{\theta}_0}$ ensures that
$\left\{h_n/(nm)\right\} {\rm Var}(\bm{\nu}'_{0n} \bm{G}_1 \bm{\nu}_{0n})
= \left\{h_n/(nm)\right\} \sigma^2_0{\rm tr}(\bm{G}_1^2  + \bm{G}_1'
\bm{G}_1)$ is bounded away from zero. By \cite{lee04}, we have
\[
\left\{h_n/(nm)\right\}^{1/2} \left\{\bm{\nu}'_{0nm} \bm{G}_1
\bm{\nu}_{0nm} - \sigma^2_0 {\rm tr}(\bm{G}_1)\right\}
\stackrel{D}{\rightarrow} N\left(0, \lim_{n \rightarrow \infty} h_n
(nm)^{-1}\sigma^4_0 \{{\rm tr}(\bm{G}_1 \bm{G}_1') + {\rm
tr}(\bm{G}_1^2)\}\right).
\]
Hence
\begin{eqnarray*}
\left\{h_n/(nm)\right\}^{1/2} \frac{\partial
\ell(\bm{\xi}_0)}{\partial \theta_1} &\stackrel{D}{\rightarrow}&
N\left(0, \lim_{n \rightarrow \infty} h_n (nm)^{-1} \left\{{\rm
tr}(\bm{G}_1 \bm{G}_1') + {\rm tr}(\bm{G}_1^2)\right\}\right)
\end{eqnarray*}
Then similarly for $\theta_k, k=2,\ldots,q$, we have
\begin{eqnarray*}
\left\{h_n/(nm)\right\}^{1/2} \frac{\partial
\ell(\bm{\xi}_0)}{\partial \theta_k} &\stackrel{D}{\rightarrow}&
N\left(0, \lim_{n \rightarrow \infty} h_n (nm)^{-1} \{{\rm
tr}(\bm{G}_k \bm{G}_k') + {\rm tr}(\bm{G}_k^2)\}\right).
\end{eqnarray*}
Thus, by the Cram\'{e}r-Wold theorem,
\[
\left\{h_n/(nm)\right\}^{1/2} \frac{\partial
\ell(\bm{\xi}_0)}{\partial\bm{\theta}}\stackrel{D}{\rightarrow}
N(\bm{0}, \bm{\Sigma}_{\bm{\theta}_0}^{-1})
\]
where
\begin{eqnarray*}
\bm{\Sigma}_{\bm{\theta}_0}^{-1} &=& \lim_{n \rightarrow \infty} - h_n (nm)^{-1} E \left\{ \frac{\partial^2 \ell(\bm{\eta}_0)}{\partial \bm{\theta} \partial \bm{\theta}'}  \right\} \\
&=& \lim_{n \rightarrow \infty} \left[\begin{array}{ccc}
h_n (nm)^{-1}{\rm tr}(\bm{G}_1^2  + \bm{G}_1' \bm{G}_1) & \ldots & h_n (nm)^{-1}{\rm tr}(\bm{G}_1 \bm{G}_q + \bm{G}_q' \bm{G}_1)\\
\vdots & \vdots & \vdots\\
h_n (nm)^{-1}{\rm tr}(\bm{G}_q \bm{G}_1 + \bm{G}_1' \bm{G}_q) &
\ldots & h_n (nm)^{-1}{\rm tr}(\bm{G}_q^2 + \bm{G}_q' \bm{G}_q)
\end{array}\right].
\end{eqnarray*}
Similarly, we can show that
\[ (nm)^{-1/2}\frac{\partial \ell(\bm{\xi}_0)}{\partial \alpha} \stackrel{D}{\rightarrow} N\left(0, \lim_{nm \rightarrow \infty} (nm)^{-1}{\rm tr}\{\bm{H}_{nm}^2+\bm{H}'_{nm}\bm{H}_{nm}\}\right) \quad {\rm and}\quad h_n(nm)^{-1}\frac{\partial{\ell^2({\bm \xi}_0)}}{\partial{{\alpha}}\partial{\bm{\theta}}} = o_p(1).\]
Thus it follows that
\begin{eqnarray*}
&& (nm/h_n)^{1/2}(\hat{\bm{\theta}}_n - \bm{\theta}_0) \\
&=& - \left\{h_n (nm)^{-1}\frac{\partial^2 \ell(\bm{\xi}_0)}{\partial \bm{\theta} \partial \bm{\theta}'}  + o_p(1) \right\}^{-1} \left\{h_n/(nm)\right\}^{1/2} \frac{\partial \ell(\bm{\xi}_0)}{\partial \bm{\theta}} -h_n^{-1/2} \left\{h_n (nm)^{-1}\frac{\partial^2 \ell(\bm{\xi}_0)}{\partial \bm{\theta} \partial \bm{\theta}'}  + o_p(1) \right\}^{-1}\\
&&\quad\quad\left\{h_n(nm)^{-1}\frac{\partial^2{\ell({\bm \xi}_0)}}{\partial{\bm{\theta}}\partial{{\alpha}}}+o_p(1) \right\}\left\{(nm)^{-1}\frac{\partial^2 \ell(\bm{\xi}_0)}{\partial \alpha^2} +o_p(1) \right\}^{-1} (nm)^{-1/2} \frac{\partial \ell(\bm{\xi}_0)}{\partial \alpha} \\
&&\stackrel{D}{\rightarrow} N(\bm{0}, \bm{\Sigma}_{\bm{\theta}_0}).
\end{eqnarray*}
and similarly,
\[(nm)^{1/2}(\hat{\alpha}^2_{nm} - \alpha^2_0) \stackrel{\mathcal{D}}{\rightarrow} N\left(0, \Sigma_{\alpha_0}\right)\]
with $\Sigma_{\alpha_0} = \lim_{nm \rightarrow \infty} nm
\left\{{\rm tr}(\bm{H}_{nm}^2+\bm{H}_{nm}'\bm{H}_{nm})
\right\}^{-1}$.

To establish the asymptotic normality of
$\hat{\bm{\beta}}_{nm}(\bm{\xi}_{nm})$ and
$\hat{\sigma}^2_{nm}(\bm{\xi}_{nm})$, we have, under (A.1)--(A.5),
\begin{eqnarray*}
&&(nm)^{1/2}(\hat{\bm{\beta}}_{nm}(\hat{\bm{\xi}}_{nm}) - \bm{\beta}_0) \\
&=& (nm)^{-1/2}\left\{(nm)^{-1} \bm{X}_{nm}' \bm{S}_{nm}'(\hat{\bm{\xi}}_{nm}) \bm{S}_{nm}(\hat{\bm{\xi}}_{nm}) \bm{X}_{nm}\right\}^{-1} \bm{X}_{nm}'\bm{S}_{nm}'(\hat{\bm{\xi}}_{nm})\bm{S}_{nm}(\hat{\bm{\xi}}_{nm})\bm{S}_{0nm}^{-1}\bm{\nu}_{0nm}\\
&=& (nm)^{-1/2}\left((nm)^{-1} \bm{X}_{nm}' \bm{S}_{0nm}' \bm{S}_{0nm} \bm{X}_{nm}\right)^{-1} \bm{X}_{nm}'\bm{S}_{0nm}'\bm{\nu}_{0nm} + o_p(1)\\
&\stackrel{D}{\rightarrow}& N\left(\bm{0}, \lim_{n\rightarrow
\infty}\left((nm)^{-1}\bm{X}_{nm}' \bm{S}_{0nm}' \bm{S}_{0nm}
\bm{X}_{nm}\right)^{-1}\right).
\end{eqnarray*}
and
\begin{eqnarray*}
(nm)^{1/2}(\hat{\sigma}^2_{nm}(\hat{\bm{\xi}}_{nm}) - \sigma^2_0)&=&
(nm)^{-1/2}\left(\bm{\nu}'_{0nm}\bm{\nu}_{0nm} - nm\sigma^2_0\right)
+ o_p(1) \stackrel{D}{\rightarrow} N(0, 2\sigma^4_0).
\end{eqnarray*}

%%%%%%%%%%%%%%%%
% Bibliography %
%%%%%%%%%%%%%%%%
\bibliographystyle{biometrics}
\bibliography{draft}

\begin{thebibliography}{}

\bibitem[\protect\bsccite{Anselin}{Anselin}{2001}]{anselin01}
Anselin, L. (2001).
\newblock Spatial econometrics.
\newblock In Baltagi, B.~H., editor, {\em A Companion to Theoretical
  Econometrics}. Blackwell Publishers Lte, Massachusetts.

\bibitem[\protect\bsccite{Baltagi}{Baltagi}{2005}]{baltagi05}
Baltagi, B.~H. (2005).
\newblock {\em Econometric Analysis of Panel Data}.
\newblock Wiley, New York, 3rd edition edition.

\bibitem[\protect\bsccite{Cressie}{Cressie}{1993}]{cressie93}
Cressie, N. (1993).
\newblock {\em Statistics for Spatial Data}.
\newblock Wiley, New York, revised edition.

\bibitem[\protect\bsccite{Horn and Johnson}{Horn and Johnson}{1985}]{horn85}
Horn, R.~A. and Johnson, C.~R. (1985).
\newblock {\em Matrix Analysis}.
\newblock Cambridge University Press.

\bibitem[\protect\bsccite{Kelejian and Prucha}{Kelejian and
  Prucha}{2001}]{kelejianp01}
Kelejian, H.~H. and Prucha, I.~R. (2001).
\newblock On the asymptotic distribution of the {M}oran {I} test statistic with
  applications.
\newblock {\em Journal of Econometrics} {\bf 104}, 217--257.

\bibitem[\protect\bsccite{Lee}{Lee}{2004}]{lee04}
Lee, L.-F. (2004).
\newblock Asymptotic distribution of quasi-maximum likelihood estimators for
  spatial autoregressive models (with supplements).
\newblock {\em Econometrica} {\bf 72}, 1899--1925.

\bibitem[\protect\bsccite{Lee and Yu}{Lee and Yu}{2010a}]{leey10a}
Lee, L.-F. and Yu, J. (2010a).
\newblock Estimation of spatial autoregressive panel data models with fixed
  effects.
\newblock {\em Journal of Econometrics} {\bf 154}, 165--185.

\bibitem[\protect\bsccite{Lee and Yu}{Lee and Yu}{2010b}]{leey10b}
Lee, L.-F. and Yu, J. (2010b).
\newblock A spatial dynamic panel data model with both time and individual
  fixed effects.
\newblock {\em Econometric Theory} {\bf 26}, 564--597.

\bibitem[\protect\bsccite{Mardia and Marshall}{Mardia and
  Marshall}{1984}]{mardiam84}
Mardia, K.~V. and Marshall, R.~J. (1984).
\newblock Maximum likelihood estimation of models for residual covariance in
  spatial regression.
\newblock {\em Biometrika} {\bf 71}, 135--146.

\bibitem[\protect\bsccite{Newey}{Newey}{1991}]{newey91}
Newey, W.~K. (1991).
\newblock Uniform convergence in probability and stochastic equicontinuity.
\newblock {\em Econometrica} {\bf 59}, 1161--1167.

\bibitem[\protect\bsccite{Robinson and Thawornkaiwong}{Robinson and
  Thawornkaiwong}{2012}]{RT12}
Robinson, P. and Thawornkaiwong, S. (2012).
\newblock Statistical inference on regression with spatial dependence.
\newblock {\em Journal of Econometrics} {\bf 167}, 521--542.

\bibitem[\protect\bsccite{White}{White}{1994}]{white94}
White, H. (1994).
\newblock {\em Estimation, Inference and Specification Analysis}.
\newblock Cambridge University Press.

\bibitem[\protect\bsccite{Yu, de~Jong and Lee}{Yu et~al.}{2008}]{yuj08}
Yu, J., de~Jong, R. and Lee, L. (2008).
\newblock Quasi-maximum likelihood estimators for spatial dynamic panel data
  with fixed effects when both $n$ and $t$ are large.
\newblock {\em Journal of Econometrics} {\bf 146}, 118--134.

\bibitem[\protect\bsccite{Zheng and Zhu}{Zheng and Zhu}{2011}]{zhengz11}
Zheng, Y. and Zhu, J. (2011).
\newblock On the asymptotics of maximum likelihood estimation for spatial
  linear models on a lattice.
\newblock {\em Sankhya - The Indian Journal of Statistics, Series A} page In
  press.

\end{thebibliography}

{\small
\newpage

\begin{table*}[!b]
 \caption{\textit{Means and standard deviations (SD) of maximum likelihood estimates (MLE) of the model parameters based on 100 simulated data. Here the lattice sizes are $4 \times 4$ and $8 \times 8$ with varying sub-lattice sizes $1 \times 1$, $2 \times 2$ and $4 \times 4$ within each cell of the lattice, and the number of time points $m$ is $2$.}}
 \label{tab:simulation1}
\begin{center}
\begin{tabular}{ccccccccc}\toprule
& & &\multicolumn{3}{c}{$4 \times 4$}&\multicolumn{3}{c}{$8 \times
8$}\\\cmidrule(r{0.5em}l{0.5em}){4-6}\cmidrule(r{0.5em}l{0.5em}){7-9}
Parameter&Truth&MLE&$1 \times 1$&$2 \times 2$&$4 \times 4$&$1 \times 1$&$2 \times 2$&$4 \times 4$\\
& & & $(n=16) $&$(n=64)$&$(n=256)$&$(n=64)$&$(n=256)
$&$(n=1024)$
\\\midrule
$\beta_0$&2.0&Mean&2.2553&1.6990&1.7385&2.1934&1.8016&2.0995\\
         &      &SD&(0.2745)&(0.0240)&(0.0188)&(0.0106)&(0.0128)&(0.0007)\\
$\beta_1$&2.0&Mean&1.8913&2.2048&2.1146&2.1913&2.1079&2.1192\\
         &      &SD&(0.3293)&(0.1026)&(0.0622)&(0.0143)&(0.0546)&(0.0254)\\
$\sigma^2$&1.0&Mean&1.4461&1.3757&1.3426&1.3532&1.3213&1.1849\\
          &      &SD&(0.3590)&(0.1890)&(0.0575)&(0.1800)&(0.0932)&(0.0262)\\
$\theta$&0.8&Mean&0.9623&0.7071&0.7300&0.6711&0.7002&0.9163\\
        &      &SD&(0.1953)&(0.0692)&(0.0528)&(0.0252)&(0.0369)&(0.0118)\\
$\alpha$&0.2&Mean&0.0225&0.0649&0.0787&0.0503&0.0820&0.1008\\
        &      &SD&(0.1389)&(0.0201)&(0.0371)&(0.0051)&(0.0167)&(0.0165)\\\midrule
\end{tabular}
\end{center}
\label{tab:simulation1}
\end{table*}

\begin{table*}[!b]
 \caption{\textit{Means and standard deviations (SD) of maximum likelihood estimates (MLE) of the model parameters based on 100 simulated data. Here the lattice size are $4 \times 4$ and $8 \times 8$ with varying sub-lattice sizes $1 \times 1$, $2 \times 2$ and $4 \times 4$ within each cell of the lattice, and the number of time points $m$ is $5$.}}
 \label{tab:simulation2}
\begin{center}
\begin{tabular}{ccccccccc}\toprule
& & &\multicolumn{3}{c}{$4 \times 4$}&\multicolumn{3}{c}{$8 \times
8$}\\\cmidrule(r{0.5em}l{0.5em}){4-6}\cmidrule(r{0.5em}l{0.5em}){7-9}
Parameter&Truth&MLE&$1 \times 1$&$2 \times 2$&$4 \times 4$&$1 \times 1$&$2 \times 2$&$4 \times 4$\\
& & & $(n=16) $&$(n=64)$&$(n=256)$&$(n=64)$&$(n=256)
$&$(n=1024)$
\\\midrule
$\beta_0$&2.0&Mean&2.2430&1.7536&2.2494&2.1473&2.1481&2.0992\\
         &      &SD&(0.0254)&(0.0058)&(0.0044)&(0.0072)&(0.0033)&(0.0028)\\
$\beta_1$&2.0&Mean&1.7381&2.1886&2.1069&2.1472&2.0865&2.0624\\
         &      &SD&(0.0616)&(0.0191)&(0.0188)&(0.0167)&(0.0220)&(0.0157)\\
$\sigma^2$&1.0&Mean&1.7640&1.2659&1.2123&1.3342&1.2075&1.1496\\
          &      &SD&(0.5033)&(0.0824)&(0.0690)&(0.1344)&(0.0309)&(0.0152)\\
$\theta$&0.8&Mean&0.7175&0.7076&0.7195&0.6907&0.7415&0.7330\\
        &      &SD&(0.0627)&(0.0246)&(0.0258)&(0.0354)&(0.0192)&(0.0171)\\
$\alpha$&0.2&Mean&0.0536&0.0774&0.0902&0.0804&0.0937&0.1131\\
        &      &SD&(0.1026)&(0.0053)&(0.0141)&(0.0188)&(0.0163)&(0.0137)\\\midrule
\end{tabular}
\end{center}
\label{tab:simulation2}
\end{table*}

\begin{table*}[!b]
 \caption{\textit{Means and standard deviations (SD) of maximum likelihood estimates (MLE) of the model parameters based on 100 simulated data. Here the lattice sizes are $4 \times 4$ and $8 \times 8$ with varying sub-lattice sizes $1 \times 1$, $2 \times 2$ and $4 \times 4$ within each cell of the lattice, and the number of time points $m$ is $10$.}}
 \label{tab:simulation3}
\begin{center}
\begin{tabular}{ccccccccc}\toprule
& & &\multicolumn{3}{c}{$4 \times 4$}&\multicolumn{3}{c}{$8 \times
8$}\\\cmidrule(r{0.5em}l{0.5em}){4-6}\cmidrule(r{0.5em}l{0.5em}){7-9}
Parameter&Truth&MLE&$1 \times 1$&$2 \times 2$&$4 \times 4$&$1 \times 1$&$2 \times 2$&$4 \times 4$\\
& & & $(n=16) $&$(n=64)$&$(n=256)$&$(n=64)$&$(n=256)
$&$(n=1024)$
\\\midrule
$\beta_0$&2.0&Mean&1.7480&2.1839&1.8039&1.8889&1.8924&1.8988\\
         &      &SD&(0.0146)&(0.0221)&(0.0161)&(0.0111)&(0.0132)&(0.0050)\\
$\beta_1$&2.0&Mean&1.8046&2.1780&1.9044&2.1789&2.1740&2.0975\\
         &      &SD&(0.0686)&(0.0319)&(0.0092)&(0.0240)&(0.0446)&(0.0105)\\
$\sigma^2$&1.0&Mean&1.3016&1.2165&1.2162&1.2812&1.1589&1.0664\\
          &      &SD&(0.1227)&(0.1708)&(0.0738)&(0.1390)&(0.1004)&(0.0670)\\
$\theta$&0.8&Mean&0.7575&0.7159&0.7158&0.7341&0.7176&0.7535\\
        &      &SD&(0.0339)&(0.0254)&(0.0212)&(0.0321)&(0.0303)&(0.0143)\\
$\alpha$&0.2&Mean&0.0654&0.0781&0.1146&0.0973&0.1331&0.1445\\
        &      &SD&(0.0675)&(0.0346)&(0.0368)&(0.0543)&(0.0329)&(0.0236)\\\midrule
\end{tabular}
\end{center}
\label{tab: simulation3}
\end{table*}
}
\end{document}